\documentclass{article}

\usepackage{fullpage}
\usepackage{url}
\usepackage{hyperref}
\usepackage[utf8]{inputenc}
\usepackage{amsthm}
\usepackage{amsmath}
\usepackage{amssymb}
\usepackage{bbm}
\usepackage{physics}
\usepackage{mathtools}
\DeclarePairedDelimiter{\ceil}{\lceil}{\rceil}
\DeclarePairedDelimiter{\floor}{\lfloor}{\rfloor}

\newtheorem{thm}{Theorem}[section]
\newtheorem{cor}[thm]{Corollary}

\newtheorem{lem}[thm]{Lemma}

\title{Minimal Gaps and Additive Energy in Real-valued Sequences}
\author{Shvo Regavim\footnote{E-mail: shvoregavim@mail.tau.ac.il}}
\date{\today}

\begin{document}
\maketitle
\begin{center}
    (Tel Aviv University)
\end{center}
\begin{abstract}
    We study the minimal gap statistic for sequences of the form $\left( \alpha x_n \right)_{n = 1}^{\infty}$ where $\left( x_n \right)_{n = 1}^{\infty}$ is a sequence of real numbers, and its connection to the additive energy of $\left( x_n \right)_{n = 1}^{\infty}$. Inspired by a recent paper of Aistleitner, El-Baz and Munsch we show conditionally on the Lindel\"{o}f Hypothesis that if the additive energy is of lowest possible order then for almost all $\alpha$, the minimal gap $\delta_{\min}^{\alpha} (N) = \min \left\{ \alpha x_m - \alpha x_n \bmod \ 1 : 1 \leq m \neq n \leq N \right\}$ is close to that of a random sequence, a result Rudnick showed for integer-valued sequences. We also show unconditional results in this direction, as well as some converse theorems about sequences with large additive energy.
\end{abstract}

\section{Introduction and Statement of Results}
Given a sequence of real numbers $(x_n)_{n \geq 1}$, our goal is to understand the minimal gap
$$\delta_{\min} \left( (x_n)_{n \geq 1}, N \right) = \min \left( \norm{x_m - x_n} : m, n \leq N, m \neq n \right)$$
Where $\norm{x} = \min_{\substack{k \in \mathbb{Z}}} \abs{x - k}$ is the distance of $x$ to the nearest integer.
In the random case, that is for $N$ independent uniformly distributed points in the unit interval, the minimal gap is almost surely of size about $\frac{1}{N^2}$. The minimal gap is always at most $\frac{1}{N}$, as that is the average size of the gap between two consecutive points.

Computing the minimal gap of a given sequence is difficult, and so we will mainly discuss the metric theory of minimal gaps. Fix a sequence $X = (x_n)_{n \geq 1}$. We will study the properties of the sequence $(\alpha x_n)_{n \geq 1}$ for $\alpha \in \mathbb{R}$. We want to understand under which conditions we can show that for almost all $\alpha$ the minimal gap
$$\delta_{\min}^{\alpha} \left( N \right) = \delta_{\min} \left( (\alpha x_n)_{n \geq 1}, N \right)$$
is of size approximately $\frac{1}{N^2}$ as in the random case. When $x_n$ are integers this problem was investigated by Rudnick \cite{Integer Minimal Gaps} and by Aistleitner, El-Baz and Munsch \cite{diff-set}. As an example of their results, Rudnick \cite{Integer Minimal Gaps} showed that if the additive energy of the sequence, a concept we define below, is $\ll N^{2 + \varepsilon}$ for all $\varepsilon > 0$, then for almost all $\alpha$ we have
$$\delta_{\min}^{\alpha} \left( N \right) \ll \frac{1}{N^{2 - \varepsilon}}$$
for all $\varepsilon > 0$. In \cite{diff-set}, Aistleitner, El-Baz and Munsch the minimal gap statistic is related to the size of the difference set, for example they show that for almost all $\alpha$, for infinitely many $N$ we have
$$\delta_{\min}^{\alpha} \left( N \right) \ll \frac{1}{C_N}$$
where $C_N = \# \left\{ x_n - x_m : 1 \leq m \neq n \leq N \right\}$ is the size of the difference set. The purpose of this paper is to develop the theory for real valued sequences, which poses several new issues.

It is relatively easy to show that under very mild spacing conditions on $x_n$ the minimal gap cannot be too small. More specifically, we have

\begin{thm}\label{lower bound}

Assume that there is some constant $c > 0$ such that $x_{n+1} - x_n \geq \frac{c \log (n)}{n}$. Then for almost all $\alpha$ we have
$$\delta_{\min}^{\alpha} \left( N \right) > \frac{1}{N^2 \log (N)^{1 + o(1)}}$$
where the implied rate of decay $o(1)$ depends on $\alpha$.

\end{thm}

To ensure that the minimal gap is small, we give the following criterion in terms of the number of solutions to a certain Diophantine inequality and a mild spacing condition on $x_n$. This follows the ideas of Rudnick and Technau \cite{Rudnick Technau}, where an analogous condition is developed for the question of the metric pair correlation of the sequence $x_n$.

\begin{thm}\label{big-sum}
Fix some $2 \geq \beta > 1$. Suppose that $x_{n+1} - x_n \geq c$ for some $c > 0$ and that for all $\varepsilon > 0$ we have
$$\sum_{1 \leq j_1, j_2 < N^{\beta - \varepsilon}} \sum_{\substack{1 \leq m_1, m_2, n_1, n_2 \leq N \\ m_1 \neq n_1, m_2 \neq n_2}} \mathbbm{1} \left( \abs{j_1 (x_{m_1} - x_{n_1}) - j_2 (x_{m_2} - x_{n_2})} < 1 \right) \leq C(\varepsilon) N^{4 - \varepsilon}$$
for some constant $C(\varepsilon)$. Then for almost all $\alpha$,
$$\delta_{\min}^{\alpha} \left( N \right) \ll \frac{1}{N^{\beta - \varepsilon}}$$
for all $\varepsilon > 0$.

\end{thm}

Our main goal in the rest of this paper is to bound this sum for various sequences. We follow the technique in \cite{2009.08184} in order to give a criterion in terms of the "additive energy" of the sequence $X$, which we define by
$$E^{*} \left( X, N, \gamma \right) = \# \left\{ (n_1, n_2, n_3, n_4) \in \left[ 1, N \right]^4 : \abs{x_{n_1} + x_{n_2} - x_{n_3} - x_{n_4}} < \gamma \right\}$$
and by convention when we do not write $\gamma$ we let it be equal to 1, that is
$$E^{*} \left( X, N \right) = E^{*} \left( X, N, 1 \right)$$
We always have $E^{*} \left( X, N \right) \geq N^2$ and by our spacing conditions it is easy to see that $E^{*} \left( X, N \right) \ll N^3$ (as $n_1, n_2, n_3$ determine $n_4$ up to a finite number of options). The additive energy measures in a sense how much additive structure our sequence has. For example, an arithmetic progression has additive energy $\gg N^3$, of highest possible order, and the sequence $x_n = 2^n$ has additive energy $\ll N^2$, of lowest possible order. In \cite{Integer Minimal Gaps} it is shown that if the additive energy of an integer-valued sequence is at most $N^{2 + o(1)}$, then the minimal gap is for almost all $\alpha$ of size $\frac{1}{N^{2 - o(1)}}$. We will, conditionally, show the corresponding result for real-valued sequences:

%Replaced beta with theta%
\begin{thm}\label{zeta v1}

Assume the Lindel\"{o}f Hypothesis. Take a sequence $X = (x_n)_{n \geq 1}$ of positive real numbers such that there exists a constant $c > 0$ such that $x_{n + 1} - x_n \geq c$, and that $x_n$ is of polynomial growth, that is there exists some constant $C$ such that $x_n \ll n^C$. Assume furthermore that for some $\theta$, for all $\varepsilon > 0$
$$E^{*} \left( X, N \right) \ll N^{\theta + \varepsilon}$$
Then for almost all $\alpha$,
$$\delta_{\min}^{\alpha} \left( N \right) < \frac{1}{N^{4 - \theta - o(1)}}$$
where the implied rate of decay $o(1)$ is dependent on $\alpha$.

\end{thm}

In particular, if the conditions of \ref{zeta v1} hold with $\theta = 2$, that is
$$E^{*} \left( X, N \right) \ll N^{2 + \varepsilon}$$
for all $\varepsilon > 0$ then for almost all $\alpha$ we have
$$\delta_{\min}^{\alpha} \left( N \right) = \frac{1}{N^{2 + o(1)}}$$
For example, if $x_n = n^{\theta}, \theta \geq 2$ then by Theorem 2 in \cite{Robert&Sargos}, $E^{*} \left( X, N \right) = N^{2 + o(1)}$ and so the theorem applies in this case. 

We will show that this result is optimal, in the sense that for all $2 \leq \theta \leq 3$ there is a sequence $x_n$ such that $E^{*} (X, N) = N^{\theta + o(1)}$ and for almost all $\alpha$ we have $\delta_{\min}^{\alpha} \left( N \right) < \frac{1}{N^{4 - \theta - o(1)}}$.

We have written this result mainly to illustrate the use of the technique. In fact, we shall prove a more general result:

\begin{thm}\label{zeta v2}

Take a sequence $X = (x_n)_{n \geq 1}$ of positive real numbers such that there exists a constant $c > 0$ such that $x_{n + 1} - x_n \geq c$. Assume furthermore that for some $\theta \geq 2$ and for all $\varepsilon > 0$
$$E^{*} \left( X, N \right) \ll N^{\theta + \varepsilon}$$
Assume that for some $\kappa \geq 2$ we have the bound
$$\frac{1}{T} \intop_{0}^{T} \abs{\zeta \left( \frac{1}{2} + i t \right)}^{\kappa} \mathrm{d} t \ll T^{m(\kappa) + \varepsilon}$$
on the $\kappa$th moment of the Riemann zeta function.
Then for almost all $\alpha$,
$$\delta_{\min}^{\alpha} \left( N \right) < \frac{1}{N^{\frac{\left( 4 - \theta \right) \left( \kappa - 2 \right)}{\kappa + 4 \cdot m(\kappa)} - o(1)}}$$
where the implied rate of decay $o(1)$ is dependent on $\alpha$.

\end{thm}

For example, suppose we knew that $\zeta \left( \frac{1}{2} + i t \right) = \mathcal{O} (t^{\mu + \varepsilon})$ for all $\varepsilon > 0$. Then, $m(\kappa) \leq \mu \kappa$ and so taking $\kappa \to \infty$ we would have a minimal gap of size at most $\frac{1}{N^{\frac{4 - \theta}{1 + 4 \mu} - \varepsilon}}$. In particular, assuming the Lindel\"{o}f Hypothesis we have $\mu = 0$, and thus we get Theorem $\ref{zeta v1}$ as a corollary.

An important special case is when the additive energy is of smallest possible order, that is $\theta = 2$. In this case directly applying Theorem \ref{zeta v2} we get

\begin{cor}\label{original zeta}

Take a sequence $X = (x_n)_{n \geq 1}$ of positive real numbers such that there exists a constant $c > 0$ such that $x_{n + 1} - x_n \geq c$. Assume furthermore that for all $\varepsilon > 0$
$$E^{*} \left( X, N \right) \ll N^{2 + \varepsilon}$$
Assume that for some $\kappa \geq 2$ we have the bound
$$\frac{1}{T} \intop_{0}^{T} \abs{\zeta \left( \frac{1}{2} + i t \right)}^{\kappa} \mathrm{d} t \ll T^{m(\kappa) + \varepsilon}$$
on the $\kappa$th moment of the Riemann zeta function.
Then for almost all $\alpha$,
$$\delta_{\min}^{\alpha} \left( N \right) < \frac{1}{N^{\frac{2 \kappa - 4}{\kappa + 4 \cdot m(\kappa)} - o(1)}}$$
where the implied rate of decay $o(1)$ is dependent on $\alpha$.

\end{cor}

It is interesting to note that the convexity bound $\zeta \left( \frac{1}{2} + i t \right) \ll t^{\frac{1}{4} + \varepsilon}$ gives the trivial bound $\frac{1}{N}$ on the minimal gap, and so any subconvexity bound on the zeta function gives a non-trivial result.

Using Heath-Brown's bound on the 12-th moment of the Riemann zeta function in \cite{12th moment} which in our terms states that $m(12) \leq 1$ we see that under the assumptions of \ref{original zeta} the minimal gap is of size at most $\frac{1}{N^{\frac{5}{4} - \varepsilon}}$. Applying Heath-Brown's bound\footnote{A tiny improvement comes from using Ivi\'{c}'s bound which is Theorem 8.3 in \cite{Ivic} where $\kappa = 178/13$ and $m(\kappa) = 16/13$, giving an exponent 1 + 31/121 = 1.2562..., and a non-trivial bound on the minimal gap for $\theta < 183/76 = 2.4079...$} in Theorem \ref{zeta v2} we get a non-trivial bound on the minimal gap for $\theta < 2.4$, that is if we know that $E^{*} (X, N) \ll N^{\theta}$ for some $\theta < 2.4$, then for some constant $\delta > 0$ the minimal gap is almost always smaller than $\frac{1}{N^{1 + \delta}}$.

Our arguments in Theorems \ref{zeta v1} and \ref{zeta v2} are in part an adaptation of those in Aistleitner, El-Baz and Munsch in \cite{2009.08184}, where these techniques are used in the setting of pair correlation.

We will also show using purely elementary arguments that certain sequences satisfy the metric property we want, specifically rapidly growing sequences and quadratic polynomials.

\begin{thm}\label{rapid}

Take a sequence $X = (x_n)_{n \geq 1}$ of positive real numbers such that for every $\varepsilon > 0$, for all $n$ large enough we have
$$x_{n + 1} > \left( 1 + \frac{1}{n^{\varepsilon}} \right) x_n$$
Then for almost all $\alpha$,
$$\delta_{\min}^{\alpha} \left( N \right) < \frac{1}{N^{2 - o(1)}}$$
where the implied rate of decay $o(1)$ is dependent on $\alpha$.

\end{thm}

Clearly every lacunary sequence satisfies this growth condition and so the conclusion of the theorem applies in this case. However, this theorem applies also to some sequences which grow subexponentially such as $x_n = e^{\frac{n}{(\log n)^{\theta}}}$ for any $\theta \in \mathbb{R}$. 

We also show this property for quadratic polynomials:

\begin{thm}\label{quadratic}

For any $a, b, c \in \mathbb{R}$ such that $a \neq 0$, look at the sequence $X = (a n^2 + b n + c)_{n \geq 1}$. For almost all $\alpha$,
$$\delta_{\min}^{\alpha} \left( N \right) < \frac{1}{N^{2 - o(1)}}$$
where the implied rate of decay $o(1)$ is dependent on $\alpha$.
\end{thm}

In the case where $a, b, c$ are integers (or rational numbers) this follows from \cite{Integer Minimal Gaps}, but the case where $a, b, c$ are general real numbers is new.

The core idea in all the upper bounds we have given so far is that the random variables $(x_n - x_m) \alpha$ for $\alpha \in \mathbb{R}$, which are the gaps, are essentially independent and therefore should behave as in the random case. In \ref{zeta v1} the additive energy is essentially a variance bound, in \ref{rapid} we use the rapid growth of $x_n$ to show the independence, and in \ref{quadratic} we use the fact that we are able to compute $x_n - x_m$ fairly explicitly.

In contrast, we shall also discuss the sequence $\sqrt{n}$, where we do not show a metric result but a result specifically about the gaps of the sequence $\sqrt{n}$ itself. Properties of the gaps of $\sqrt{n}$ have been studied in \cite{Elkies McMullen} where their gap distribution is computed and is surprisingly shown to not be Poissonian (as in the random case), and in \cite{Correlation sqrt} where in contrast their pair correlation is shown to be Poissonian. In the result on pair correlation, as in our case, we have to throw away the squares as they trivially cause the pair correlation to fail or the minimal gap to be $0$. We will compute the minimal gap of the sequence:

\begin{thm}\label{square root}

The minimal gap of the sequence $\left( \sqrt{n} \right)$ is asymptotically equal to $\frac{1}{2 N^{\frac{3}{2}}}$, that is
$$\delta_{\min} \left( \left( \sqrt{n} \right)_{k^2 \neq n}, N \right) \sim \min \left( \norm{\sqrt{m} - \sqrt{n}} : m, n \leq N, m \neq n, \sqrt{m}, \sqrt{n} \not\in \mathbb{Z} \right) \sim \frac{1}{2 N^{\frac{3}{2}}}$$

\end{thm}

Notice that this is in stark contrast to the random case.

In the last section we discuss sequences with large additive energy. We will try to understand how true the converse to \ref{zeta v2} is: given a sequence with large additive energy, what can we say about its minimal gap? It is easy to generate a sequence with large additive energy, but small metric minimal gap, simply by interlacing the terms of a sequence with large additive energy and a sequence with small metric minimal gap. Therefore, the best one can hope for given a sequence with large additive energy, is to find a relatively dense subsequence with large minimal gap. In the case where the sequence is integer-valued we establish the following:

\begin{thm}\label{Large Energy}

Let $x_n$ be a monotone increasing sequence of positive integers such that for some $r > 0$
$$E_N := \# \left\{ (a, b, c, d) \in [1, N]^4 \ \big\rvert \ x_a + x_b = x_c + x_d \right\} \geq r N^3$$
Then there exist increasing sequences of positive integers $n_k, a_k \to \infty$ such that $n_k$ is of density bounded from below in $\left\{ 1, 2, \dots, a_N \right\}$ and for almost all $\alpha$ we have
$$\delta_{\min}^{\alpha} (x_{n_k}, a_N) \gg \frac{1}{a_N \log (a_N)^{1 + \varepsilon}}$$

\end{thm}

The key is the Balog-Szemer\'{e}di-Gowers theorem, which essentially says that a set with large additive energy contains a large subset with a relatively small difference set, and we shall show that an integer sequence with a small difference set has a large minimal gap. For the case of real-valued sequences, the answer is more subtle and depends on the covering number of the difference set with respect to different scales. Define the covering number of a set $A \in \mathbb{R}$ with respect to a scale $\gamma$ to be
$$\mathrm{cov} (A, \gamma) = \inf_{B \subset \mathbb{R}} \left\{ \abs{B} \ \big\rvert \ \forall a \in A \ \exists b \in B: \abs{a - b} < \gamma \right\}$$
that is, the smallest number of balls with a radius of $\gamma$ that are needed to cover $A$. A natural question is what are the conditions on $\gamma$ required in order to deduce that a sequence whose difference set has a small covering number on scale $\gamma$ has a large minimal gap. In the case where $\gamma$ is small we prove

\begin{thm}\label{Fine Scale}

For any sequence of real numbers $x_n$ such that $x_{n + 1} - x_n \geq \frac{c \log n}{n}$ for some constant $c > 0$, if for some $\varepsilon_0 > 0$ we have
$$\mathrm{cov} \left( X_N - X_N, \gamma \right) \ll N$$
where $\gamma = \frac{1}{N \log(N)^{1 + \varepsilon_0}}$, we have that for every $\varepsilon > 0$ and for almost all $\alpha \in \mathbb{R}$ that 
$$\delta_{\min}^{\alpha} (x_n) \gg \frac{1}{N \log(N)^{1 + \varepsilon}}$$

\end{thm}

From this we will deduce the following corollary, which is an analogue of Theorem \ref{Large Energy} in the real-valued case:

\begin{cor}\label{Corr}

Let $x_n$ be a sequence of real numbers such that $x_{n + 1} - x_n \geq \frac{c \log n}{n}$ for some constant $c > 0$, and for some constant $r > 0$ we have
$$E \left( X, N, \gamma \right) \geq r N^3$$
Then there exist increasing sequences of positive integers $n_k, a_k \to \infty$ such that $n_k$ is of density bounded from below in $\left\{ 1, \dots, a_N \right\}$ and for almost all $\alpha$ we have
$$\delta_{\min}^{\alpha} (x_{n_k}, a_N) \gg \frac{1}{a_N \log(a_N)^{1 + \varepsilon}}$$

\end{cor}

If $\gamma$ is slightly larger, we show:

\begin{thm}\label{Coarse Scale}

There exists a sequence $\varepsilon_n \downarrow 0$ such that the sequence $x_n = n - \frac{1}{n^{1 - \varepsilon_n}}$ satisfies for all $\varepsilon > 0$ and for almost all $\alpha \in \mathbb{R}$ that
$$\delta_{\min}^{\alpha} (N) \ll \frac{1}{N^{2 - \varepsilon}}$$
where the implied constant depends on $\alpha$ and $\varepsilon$.

Notice that this sequence satisfies that
$$\mathrm{cov} \left( X_N - X_N, \gamma \right) \ll N$$
for $\gamma = \frac{1}{N^{1 - \varepsilon}}$

\end{thm}

As we will briefly discuss when we prove this theorem, this is a phenomenon that occurs only in the real valued case, and is not visible when dealing with integer valued sequences. For converse results on the absence of almost sure Poisson pair correlation if the additive energy is large, see \cite{Converse Aist, Converse Bloom}.

\section{Preliminaries}
We will use the notation $\mathbbm{1}_{[- \frac{1}{M}, \frac{1}{M}]}$ for the indicator function of the interval $[-1/M, 1/M]$ extended with period 1. Let us look at the function
$$X_{N, M} (\alpha) = \sum_{1 \leq m \neq n \leq N} \mathbbm{1}_{[- \frac{1}{M}, \frac{1}{M}]} \left( \alpha x_m - \alpha x_n \right)$$
By definition, $X_{N, M} (\alpha)$ is positive if and only if $\delta_{\min}^{\alpha}(x_n) \leq \frac{1}{M}$, and so our question is for what size of $M$ can we ensure that $X_{N, M} (\alpha)$ is almost always positive. The general idea is to find the expectation of $X_{N, M}$, bound its variance and then finish by a standard argument using Chebyshev's inequality and Borel-Cantelli. The variance will turn out to be essentially the number of solutions to a certain Diophantine inequality. We will first bound this in the general case using the Riemann zeta function, and in two specific cases (quadratic polynomials and rapidly growing sequences) we will bound the sum directly in an elementary way.

In the integer-valued case, we can simply integrate with respect to the Lebesgue integral on the interval $[0, 1]$, and use the fact that ${(\alpha + 1) x_n} = {\alpha x_n}$. In the real-valued case we cannot do this as the Lebesgue measure is not a probability measure on $\mathbb{R}$, and so we must integrate with respect to a different measure. A good choice is the measure $\mu$ which is defined by
$$\mathrm{d} \mu (x) = a(x) \mathrm{d} x := \frac{\sin^{2} \left( \pi x \right)}{\left( \pi x \right)^2} \mathrm{d} x$$
Given a function $F$, our definition of the Fourier transform will be
$$\hat{F}(\xi) = \int_\mathbb{R} e^{-i \xi t} F(t) \mathrm{d} t$$
The Fourier transform of the function $a(x) = \frac{\sin \left( \pi x \right)^2}{\pi^2 x^2}$ is
$$\hat{a}(y)  = \max \left( \frac{2 \pi - \abs{y}}{2 \pi}, 0 \right)$$
which is compactly supported on the interval $(-2 \pi, 2 \pi)$, and uniformly bounded by $1$. Notice that we have normalized $a$ so that
$$\int_\mathbb{R} a(x) \mathrm{d} x = 1$$

To make the analysis cleaner it will be useful to replace $\mathbbm{1}_{[- \frac{1}{M}, \frac{1}{M}]}$ with a smoother function. A nice choice are the Selberg polynomials $f_{M-1, 1, M}^{\pm}$, which are trigonometric polynomials of degree at most $M-1$ that satisfy
$$f_{M-1, 1, M}^{-} (x) \leq \mathbbm{1}_{[- \frac{1}{M}, \frac{1}{M}]} (x) \leq f_{M-1, 1, M}^{+} (x)$$
for all $x$, and also that
$$f_{M-1, 1, M}^{\pm} (x) = \sum_{\abs{j} < M} c_{j}^{\pm} e^{2 \pi i j x}$$
where
$$\abs{c_{j}^{\pm}} \leq \min \left( \frac{2}{M}, \frac{1}{\pi \abs{j}} \right) + \frac{1}{M} \leq \frac{3}{M}$$
and also
$$c_{0}^{-} = \intop_{0}^{1} f_{M - 1, 1, M}^{-} (x) \mathrm{d} x = \frac{1}{M}, \ c_{0}^{+} = \intop_{0}^{1} f_{M - 1, 1, M}^{+} (x) \mathrm{d} x = \frac{3}{M}$$
Now we can define
$$X_{N, M}^{\pm} (\alpha) = \sum_{1 \leq m \neq n \leq N} f_{M - 1, 1, M}^{\pm} \left( \alpha x_m - \alpha x_n \right)$$
Notice that
$$X_{N, M}^{-} (\alpha) \leq X_{N, M} (\alpha) \leq X_{N, M}^{+} (\alpha)$$
Therefore bounding $X_{N, M}^{+}$ from above corresponds to bounding the minimal gap from below, and bounding $X_{N, M}^{-}$ from below corresponds to bounding the minimal gap from above.
For clarity, we will write $f^{\pm} := f_{M - 1, 1, M}^{\pm}$.

\section{Minimal Gaps are not too Small}\label{Lower Bound}
In this section we shall prove Theorem \ref{lower bound}. We shall do this by computing the expectation of $X_{N, M}^{+}$ and showing that it is small. Notice that
$$\abs{\mathbbm{E} \left[ X_{N, M}^{+}, \mu \right] - \frac{3 N(N - 1)}{M}} = \abs{\sum_{1 \leq m \neq n \leq N} \int_{\mathbb{R}} \left( f^{+} \left( \alpha x_m - \alpha x_n \right) - \frac{3}{M} \right) a(\alpha) \mathrm{d} \alpha} =$$
$$= \abs{\sum_{0 < \abs{j} < M} c_{j}^{+} \sum_{1 \leq m \neq n \leq N} \int_{\mathbb{R}} e^{2 \pi i j \left( \alpha x_m - \alpha x_n \right)} a(\alpha) \mathrm{d} \alpha} \leq$$
$$\leq \frac{3}{M} \sum_{0 < \abs{j} < M} \sum_{1 \leq m \neq n \leq N} \hat{a} (2 \pi j (x_n - x_m)) \leq$$
$$\leq \frac{3}{M} \sum_{1 \leq m \neq n \leq N} \sum_{0 < \abs{j} < M} \mathbbm{1} \left( \abs{j (x_n - x_m)} < 1 \right)$$
Assume without loss of generality that $n > m$. By our spacing conditions, for a pair $j, n, m$ to contribute to the sum we must have
$$1 > \abs{j (x_n - x_m)} > \abs{\frac{j (n - m) \log (n)}{n}}$$
and so we must have
$$\frac{n}{\log (n)} > j (n - m)$$
Fixing $n$, the number on the right is the number of integer points under the hyperbola $x y = \frac{n}{\log (n)}$ and this is easily seen to be $\ll n$. Summing over $n$ from $1$ to $N$ we get that the error term contributes $\ll \frac{N^2}{M}$ and so
$$\mathbbm{E} \left[ X_{N, M}^{+}, \mu \right] \ll \frac{N^2}{M}$$

Fix some constant $\varepsilon > 0$, and for every $\gamma \in \mathbb{R}, N \in \mathbb{N}$ define
$$A (N, \gamma) = \left\{ \alpha \in \mathbb{R} : a(\alpha) \geq \gamma, \delta_{\min}^{\alpha} (N) \leq \frac{1}{N^2 \log (N)^{1 + \varepsilon}} \right\}$$

Using our bound on the expectation we get
$$\gamma \cdot \abs{A (N, \gamma)} \leq \mathbbm{E} \left[ X_{N, N^2 \log (N)^{1 + \varepsilon}}^{+}, \mu \right] \ll \frac{1}{\log (N)^{1 + \varepsilon}}$$
and so
$$\abs{A (N, \gamma)} \ll \frac{1}{\gamma \log (N)^{1 + \varepsilon}}$$
Using Borel-Cantelli and the fact that the series $\sum_{n = 1}^{\infty} \frac{1}{n^{1 + \varepsilon}}$ converges, we see that the set of $\alpha$ that are in $A (2^N, \gamma)$ for infinitely many $N$ is of measure 0. Take $\alpha$ such that for all $N$ large enough, $\alpha \not \in A (2^N, \gamma)$. For such an $\alpha$, for any large enough $n \in \mathbb{N}$, taking $N$ to be the unique positive integer such that
$$2^{N - 1} < n \leq 2^N$$
we see that
$$\delta_{\min}^{\alpha} (n) \geq \delta_{\min}^{\alpha} (2^N) \geq \frac{1}{2^{4N} \log (2^N)^{1 + \varepsilon}} \geq \frac{1}{5 n^2 \log (n)^{1 + \varepsilon}}$$
Therefore, for any $\gamma > 0$, the set of $\alpha$ that satisfy
$$a (\alpha) \geq \gamma$$
and also do not satisfy
$$\delta_{\min}^{\alpha} (n) \geq \frac{1}{5 n^2 \log (n)^{1 + \varepsilon}}$$
is of measure 0. As the set of zeroes of $a (\alpha)$ is countable, and therefore of measure 0, we see that for almost every $\alpha$, for all $n$ large enough we have
$$\delta_{\min}^{\alpha} (n) \geq \frac{1}{5 n^2 \log (n)^{1 + \varepsilon}}$$
for every $\varepsilon > 0$, as required.

\section{Minimal Gaps can be Small}
In this section, we will show how the variance of $X_{N, M}^{\pm}$ can be bounded in terms of the sum appearing in Theorem \ref{big-sum}. From here on we assume that for some constant $c > 0$ we have $x_{n + 1} - x_n \geq c$. In fact, rescaling $x_n$ we can assume that $x_{n + 1} - x_n \geq 1$, as this obviously does not affect metric properties. In this case we have
$$\mathbbm{E} \left[ X_{N, M}^{+}, \mu \right] = \frac{3 N(N - 1)}{M}$$
$$\mathbbm{E} \left[ X_{N, M}^{-}, \mu \right] = \frac{N(N - 1)}{M}$$
as looking at our computations in the previous section we see that $\abs{2 \pi j (x_n - x_m)} \geq 2 \pi$, and so $\hat{a} \left( 2 \pi j (x_n - x_m) \right) = 0$ which implies this result.

Now we will bound the variance of $X_{N, M}^{-}$. By definition, we have
$$\mathrm{Var} \left( X_{N, M}^{-}, \mu \right) = \mathrm{Var} \left( X_{N, M}^{-} - \frac{N(N -1)}{M}, \mu \right) = \mathbbm{E} \left[ \left( X_{N, M}^{-} - \frac{N(N - 1)}{M} \right)^2, \mu \right] =$$
$$= \sum_{0 < \abs{j_1}, \abs{j_2} < M} \sum_{\substack{m_1 \neq n_1, m_2 \neq n_2 \\ 1 \leq m_1, n_1, m_2, n_2 \leq N}} c_{j_1}^{-} c_{j_2}^{-} \int_{\mathbb{R}} e^{2 \pi i \alpha \left( j_1 \left( x_{m_1} - x_{n_1} \right) - j_2 \left( x_{m_2} - x_{n_2} \right) \right)} a(\alpha) \mathrm{d} \alpha \ll$$
$$\ll \frac{1}{M^2} \sum_{0 < j_1, j_2 < M} \sum_{\substack{m_1 < n_1, m_2 < n_2 \\ 1 \leq m_1, n_1, m_2, n_2 \leq N}} \mathbbm{1} \left( \abs{j_1 \left( x_{n_1} - x_{m_1} \right) - j_2 \left( x_{n_2} - x_{m_2} \right)} < 1 \right)$$
For technical reasons, it will be more convenient to use the Cauchy-Schwarz inequality to separate $j_1, j_2$ into dyadic intervals. Writing $U = \ceil{\log_2 (M)}$, so that $2^U > M$ we get
$$\mathrm{Var} \left( X_{N, M}^{-}, \mu \right) \ll \int_{\mathbb{R}} \left( \sum_{u = 1}^{U} \sum_{1 \leq m \neq n \leq N} \sum_{2^{u - 1} \leq j < 2^u} c_{j}^{-} e^{2 \pi i j \alpha (x_m - x_n)} \right)^2 \mathrm{d} \mu <$$
$$< \int_{\mathbb{R}} \left( \sum_{k = 1}^{U} 1 \right) \sum_{u = 1}^{U} \abs{\sum_{2^{u - 1} \leq \abs{j} < 2^u} \sum_{1 \leq m \neq n \leq N} c_{j}^{-} e^{2 \pi i j \alpha (x_m - x_n)}}^2 a(\alpha) \mathrm{d} \alpha \ll$$
$$\ll \frac{\log (N)}{M^2} \sum_{u = 1}^{U} \sum_{2^{u - 1} \leq j_1, j_2 < 2^u} \sum_{\substack{m_1 < n_1, m_2 < n_2 \\ 1 \leq m_1, n_1, m_2, n_2 \leq N}} \mathbbm{1} \left( \abs{j_1 \left( x_{n_1} - x_{m_1} \right) - j_2 \left( x_{n_2} - x_{m_2} \right)} < 1 \right)$$
Recall that we have assumed that this sum is $\ll N^{4 - \varepsilon}$ if $M = N^{\beta - \varepsilon}$, and therefore
$$\mathrm{Var} \left( X_{N, M}^{-}, \mu \right) \ll \frac{N^{4 - \varepsilon}}{M^2}$$
In the next section we will show how a small variance implies a relatively small minimal gap which will complete the proof of Theorem \ref{big-sum}, and then we shall bound the number of solutions to this Diophantine inequality in various cases.

\section{Variance to Minimal Gap}
This argument is standard but is included for the sake of completeness. In general, suppose we are given a sequence of random variables $Z_n$ such that $\mathbbm{E} \left[ Z_n \right] > 0$ and we wish to show that almost always, all but finitely many $Z_n$ are positive. Notice that by Chebyshev's inequality we have

$$\mathbb{P} \left( Z_n \leq 0 \right) \leq \mathbb{P} \left( \abs{Z_n - \mathbbm{E} \left[ Z_n \right]} \geq \mathbbm{E} \left[ Z_n \right] \right) \leq \frac{\mathrm{Var} (Z_n)}{\mathbbm{E} \left[ Z_n \right]^2}$$
and so by Borel-Cantelli it is enough that
$$\sum_{n} \frac{\mathrm{Var} (Z_n)}{\mathbbm{E} \left[ Z_n \right]^2} < \infty$$
In our case, taking $M = N^{\beta - \varepsilon}$ and using what we know about the variance and the expectation we have
$$\frac{\mathrm{Var} \left( X_{N, M}^{-}, \mu \right)}{\mathbbm{E} \left[ X_{N, M}^{-}, \mu \right]^2} \ll \frac{1}{N^{\varepsilon}}$$
Now taking the random variables
$$Z_n = X_{2^n, 2^{n (\beta - \varepsilon)}}$$
we see that
$$\sum_{n} \frac{\mathrm{Var} (Z_n)}{\mathbbm{E} \left[ Z_n \right]^2} \ll \sum_{n} \frac{1}{2^{\varepsilon n}} < \infty$$
and so for almost all $\alpha$, for all $n$ large enough we have
$$\delta_{\min}^{\alpha} \left( 2^n \right) < \frac{1}{2^{n (\beta - \varepsilon)}}$$
Now, taking any $N$ large enough, and choosing $n$ to be the unique integer such that $2^n \leq N < 2^{n + 1}$
we have
$$\delta_{\min}^{\alpha} (N) \ll \frac{1}{N^{\beta - \varepsilon}}$$
for almost all $\alpha$. Now, taking $\varepsilon = \frac{1}{n}$ and taking a countable intersection of sets of full measure we see that for almost all $\alpha$,
$$\delta_{\min}^{\alpha} (N) < \frac{1}{N^{\beta - o(1)}}$$
as required.

\section{Bounding the Riemann zeta function}
In the following sections, we will prove Theorems \ref{zeta v1} and \ref{zeta v2}, following the method in Aistleitner, El-Baz and Munsch \cite{2009.08184}. This is slightly technical, so let us explain the general idea:

We are trying to bound the amount of pairs of numbers of the form $j_1 \left( x_{n_1} - x_{m_1} \right)$ and $j_2 \left( x_{n_2} - x_{m_2} \right)$ which are close. This is the same as asking when $\frac{j_1 \left( x_{n_1} - x_{m_1} \right)}{j_2 \left( X_{n_2} - x_{m_2} \right)}$ is close to 1. In general, given a sequence of numbers $y_1, y_2, \dots$, in order to find how many numbers are close to $1$, we can look at the Dirichlet polynomial $\sum y_{n}^{i t}$, and integrate this multiplied by a function $\Phi$ whose Fourier transform is supported in a small interval around 0. In our case, the multiplicative structure of the sequence $y_n$ allows us to factor the Dirichlet polynomial as
$$\abs{\sum j^{i t}}^2 \abs{\sum \left( x_n - x_m \right)^{i t}}^2$$
The first sum can be related to the Riemann zeta function, and then we can bound it using known bounds on the Riemann zeta function. As for the second sum, when we multiply it by $\Phi$ and integrate it we get the number of pairs of numbers of the form $x_{n_1} - x_{m_1}$ and $x_{n_2} - x_{m_2}$ which are close to each other, which is essentially the additive energy.
\\ \\
In this section we will prove \ref{zeta v1} as an introduction to the method, and leave \ref{zeta v2} for the next section, as there are some extra technicalities.
%We shall now prove Theorems \ref{zeta v1} and \ref{zeta v2}, following the method in Aistleitner, El-Baz and Munsch \cite{2009.08184}.%
\\
To simplify notation, let $\left\{ z_1, \dots, z_{\frac{N (N - 1)}{2}} \right\}$ be the multiset (that is, multiplicity is allowed) of differences $\left\{ x_n - x_m : 1 \leq m < n \leq N \right\}$. Then, from the computation in the previous sections we see that it is enough to show that the sum
$$\sum_{2^{u - 1} \leq j_1, j_2 < 2^u} \sum_{1 \leq m, n \leq N} \mathbbm{1} (\abs{j_1 z_m - j_2 z_n} < 1)$$
is bounded by $N^{4 - \varepsilon}$ for $2^u \leq N^{4 - \theta - \varepsilon}$. Notice that $\abs{j_1 z_m - j_2 z_n} < 1$ is equivalent to $\abs{\frac{j_1}{j_2} - \frac{z_n}{z_m}} < \frac{1}{j_2 z_m}$ and so by the triangle inequality $\frac{z_n}{z_m} < \frac{j_1 + 1}{j_2} \leq 2$. In the same way we get $\frac{1}{2} < \frac{z_n}{z_m}$.

Since we assumed that $x_n$ is of polynomial growth, we can divide the interval $\left[ 1, N^{C} \right]$ into approximately $\log N$ dyadic intervals of the form $\left[ 2 N^{\beta}, 4 N^{\beta} \right]$ where $\beta \in \left[ 0, C \right]$, and split the sum according to the dyadic interval that $z_n$ lies in. As we have just shown above, if $z_n \in \left[ 2 N^{\beta}, 4 N^{\beta} \right]$ then $z_m \in \left[ N^{\beta}, 8 N^{\beta} \right]$, and so we see that it is enough to show that the sum
$$\sum_{2^{u - 1} \leq j_1, j_2 < 2^u} \sum_{z_m, z_n \in [N^{\beta}, 8 N^{\beta}]} \mathbbm{1} (\abs{j_1 z_m - j_2 z_n} < 1)$$
is bounded by $N^{4 - \varepsilon}$ for $\beta \in \left[ 0, C \right]$.

%We don't really need to assume that beta is not small?%
%Using this, we can bound the contribution of small $N^{\beta}$: notice that fixing $j_1, z_m, z_n$ there is at most 1 option for $j_2$, and thus there are at most $2^u N^{1 + \beta + \varepsilon} \ll N^{3 + \beta + \varepsilon}$ solutions. Thus, the contribution of $\beta < \frac{3}{4}$ is at most $N^{3 + \frac{3}{4} + \varepsilon}$, which is small, and so from here on we will assume $\beta \geq \frac{3}{4}$.%

Let $T = 2^{u - 1} N^{\beta}$. The inequality $\abs{j_1 z_m - j_2 z_n} < 1$ implies
$$\abs{\frac{j_1 z_m}{j_2 z_n} - 1} < \frac{1}{j_2 z_n} \leq \frac{1}{2^{u - 1} N^{\beta}} = \frac{1}{T}$$
For $x$ small enough we have $\abs{\log ( 1 + x )} \leq 2 \abs{x}$ and so we just need to bound the number of solutions to
$$\abs{\log \left( \frac{j_1 z_m}{j_2 z_n} \right)} < \frac{2}{T}$$
Define the Dirichlet polynomials
$$D(t) = \sum_{2^{u - 1} \leq j < 2^u} j^{i t}$$
$$P(t) = \sum_{z_n \in \left[ N^{\beta}, 8 N^{\beta} \right]} z_{n}^{i t}$$
We see that
$$\abs{D(t)}^2 \abs{P(t)}^2 = \sum_{2^{u - 1} \leq j_1, j_2 < 2^u} \sum_{z_m, z_n \in \left[ N^{\beta}, 8 N^{\beta} \right]} \left( \frac{j_1 z_m}{j_2 z_n} \right)^{i t}$$
Take $\Phi$ to be a Schwarz function whose Fourier transform $\hat{\Phi}$ is non-negative, compactly supported in $[-4, 4]$ and bounded there by $2$, and for $x \in [-2, 2]$ we have $\hat{\Phi} (x) \geq 1$. Then, notice that
$$\sum_{2^{u - 1} \leq j_1, j_2 < 2^u} \sum_{z_m, z_n \in [N^{\beta}, 8 N^{\beta}]} \mathbbm{1} (\abs{j_1 z_m - j_2 z_n} < 1) \leq \int_{\mathbb{R}} \abs{D(t T)}^2 \abs{P(t T)}^2 \Phi(t) \mathrm{d} t$$
Notice that we have the pointwise bounds
$$\abs{D(t)}^2 \leq D(0)^2 \leq 2^{2 u}$$
$$\abs{P(t)}^2 \leq P(0)^2 = \left( \sum_{z_n \in \left[ N^{\beta}, 8 N^{\beta} \right]} 1 \right)^2$$
Obviously we have
$$\sum_{z_n \in \left[ N^{\beta}, 8 N^{\beta} \right]} 1 \leq N^2$$
but by Cauchy-Schwarz we see that
$$\sum_{z_n \in \left[ N^{\beta}, 8 N^{\beta} \right]} 1 \leq \sum_{k \in \left[ N^{\beta}, 8 N^{\beta} \right]} \left( \sum_{n} \mathbbm{1} \left( z_n \in [k, k + 1) \right) \right) \ll \left( N^{\beta} E^{*} (X, N) \right)^{\frac{1}{2}}$$
which gives a bound of
$$\abs{P(t)}^2 \ll \min \left( N^4, N^{\beta + \theta + \varepsilon} \right)$$
Using this, we can bound the contribution of $\abs{t} \leq \frac{1}{T}$ by
$$\frac{N^{\beta + \theta + \varepsilon} \cdot 2^{2 u}}{T} \ll 2^{u} N^{\theta} \ll N^{4 - \varepsilon}$$
which is small enough. For $\abs{t} \geq \frac{1}{T}$ we will connect $\abs{D(t)}^2$ to the Riemann zeta function via a convolution formula. This lemma is a special case of Lemma 5.3 in \cite{Gal-Sums}.

\begin{lem}\label{lem1}

Let $F$ be a holomorphic function in the strip $y = \Im(z) \in \left[ -\frac{3}{2}, 0 \right]$ such that
$$\sup_{-\frac{3}{2} \leq y \leq 0} \abs{F(x + i y)} \ll \frac{1}{x^2 + 1}$$
Then for all $s = \frac{1}{2} + i t, t \neq 0$ we have
$$\sum_{k, \ell \geq 1} \frac{\hat{F} (\log k \ell)}{k^s \ell^{\bar{s}}} = \int_{\mathbb{R}} \zeta(s + i u) \overline{\zeta(s - i u)} F(u) \mathrm{d} u + 2 \pi \zeta(1 - 2 i t) F(i s - i) + 2 \pi \zeta(1 + 2 i t) F(i \bar{s} - i)$$

\end{lem}

We will sketch the proof: Integrate the function $w \mapsto \zeta(s + w) \zeta(\bar{s} + w) F(-i w)$ on the vertical line $\Re(w) = 0$. Shift the contour to $\Re(w) = \frac{3}{2}$, and the residues that are collected correspond to the terms $2 \pi \zeta(1 - 2 i t) F(i s - i) + 2 \pi \zeta(1 + 2 i t) F(i \bar{s} - i)$. On this line, the series for $\zeta$ converges absolutely, and we can interchange the order of summation and integration. Then, shifting termwise each integral back to $\Re(w) = 0$ we get exactly our identity.

In our case, we will take $F(x) = \frac{\sin^2 \left( K x \right)}{\pi K x^2}$ where $K = \left( 1 + \varepsilon \right) \log \left( 2^u \right)$. Notice that $\hat{F} (y) = \min \left( 0, 1 - \frac{\abs{y}}{2 K} \right)$, and also that $F \left( t \pm \frac{i}{2} \right) \ll \frac{2^{u \left( 1 + \varepsilon \right)}}{t^2 + 1}$. Then Lemma \ref{lem1} gives
$$\sum_{j_1, j_2 \geq 1} \frac{\hat{F} \left( \log j_1 j_2 \right)}{(j_1 j_2)^{\frac{1}{2}}} \left( \frac{j_1}{j_2} \right)^{i t} =$$
$$= \int_{\mathbb{R}} \zeta \left( \frac{1}{2} + i t + i u \right) \zeta \left( \frac{1}{2} - i t + i u \right) F(u) \mathrm{d} u + 2 \pi \zeta \left( 1 - 2 i t \right) F \left( - t - \frac{i}{2} \right) + 2 \pi \zeta \left( 1 + 2 i t \right) F \left( t - \frac{i}{2} \right)$$
Since $\hat{\Phi}$ is non-negative, increasing the coefficients of $\abs{D(t)}^2$ can only increase our integral, and so we can replace $\abs{D(t)}^2$ with
$$2^u \sum_{j_1, j_2 \geq 1} \frac{\hat{F} \left( \log j_1 j_2 \right)}{(j_1 j_2)^{\frac{1}{2}}} \left( \frac{j_1}{j_2} \right)^{i t}$$
Now using the convolution identity we get
$$\int_{\abs{t} \geq \frac{1}{T}} \abs{D(t T)}^2 \abs{P(t T)}^2 \Phi (t) \mathrm{d} t \ll 2^u \left( \mathrm{I_1} + \mathrm{I_2} \right)$$
where
$$\mathrm{I_1} = \int_{t \geq \frac{1}{T}} \int_{\mathbb{R}} \zeta \left( \frac{1}{2} + i t T + i u \right) \zeta \left( \frac{1}{2} - i t T + i u \right) F(u) \ \mathrm{d} u \abs{P(t T)}^2 \Phi(t) \ \mathrm{d} t$$
$$\mathrm{I_2} = \int_{t \geq \frac{1}{T}} \zeta \left( 1 + 2 i t T \right) F \left( t T - \frac{i}{2} \right) \abs{P(t T)}^2 \Phi(t) \ \mathrm{d} t$$
Using the easy bound $\zeta \left( 1 + i t \right) \ll \log(t + 1)$ for $t \geq 1$ we have
$$\mathrm{I_2} \ll \frac{2^{2u \left( 1 + \varepsilon \right)} N^{\beta}}{T} \int_{t \geq 1} \frac{\log (t + 1)}{t^2 + 1} \Phi \left( \frac{t}{T} \right) \mathrm{d} t \ll 2^{u \left( 1 + \varepsilon \right)}$$

Let us just mention that later we bound the integral of $\abs{P(t T)}^2 \Phi (t)$ in a more efficient way and so we can get a much better bound on $I_2$, but this isn't necessary for our purposes. Now we turn to bounding $\mathrm{I_1}$. For simplicity, let us assume the Lindel\"{o}f Hypothesis, that is $\abs{\zeta \left( \frac{1}{2} + i t \right)} \ll t^{\varepsilon}$ for all $\varepsilon > 0$. Then,
$$\abs{\zeta \left( \frac{1}{2} + i t T + i u \right) \zeta \left( \frac{1}{2} - i t T + i u \right)} \ll \left( T^{\frac{\varepsilon}{2}} \abs{t^{\frac{\varepsilon}{2}}} + u^{\frac{\varepsilon}{2}} \right)^{2} \ll T^{\varepsilon} t^{\varepsilon} + u^{\varepsilon} \ll N^{\varepsilon} \left( t^{\varepsilon} + u^{\varepsilon} \right)$$
and thus
$$\mathrm{I_1} \ll N^{\varepsilon} \left( \mathrm{Int_a} + \mathrm{Int_b} \right)$$
where
$$\mathrm{Int_a} = \int_{t \geq \frac{1}{T}} \int_{\mathbb{R}} u^{\varepsilon} F(u) \ \mathrm{d} u \abs{P(t T)}^2 \Phi (t) \ \mathrm{d} t$$
$$\mathrm{Int_b} = \int_{t \geq \frac{1}{T}} \int_{\mathbb{R}} F(u) \mathrm{d} u \ t^{\varepsilon} \abs{P(t T)}^2 \Phi(t) \ \mathrm{d} t$$
As the integral of $u^{\varepsilon} F(u) \ll \frac{u^{\varepsilon}}{u^2 + 1}$ is bounded uniformly for $0 \leq \varepsilon \leq 1 - c < 1$, we just need to bound the integral over $t$ of $\abs{P(t T)}^2 \Phi(t)$ and $t^{\varepsilon} \abs{P(t T)}^2 \Phi(t)$. By the rapid decay of $\Phi$, the contribution of $t > N^{\varepsilon}$ is negligible, and thus both $\mathrm{Int_a}$ and $\mathrm{Int_b}$ are bounded by
$$N^{\varepsilon} \int_{\mathbb{R}} \abs{P(t T)}^2 \Phi(t) \ \mathrm{d} t \leq 2 N^{\varepsilon} \sum_{z_m, z_n \in \left[ N^{\beta}, 8 N^{\beta} \right]} \mathbbm{1} \left( \abs{\log \frac{z_m}{z_n}} < \frac{4}{T} \right) \ll$$
$$\ll N^{\varepsilon} \sum_{z_m, z_n \in \left[ N^{\beta}, 8 N^{\beta} \right]} \mathbbm{1} \left( \abs{\frac{z_m}{z_n} - 1} < \frac{8}{2^{u - 1} N^{\beta}} \right) \leq N^{\varepsilon} \sum_{z_m, z_n \in \left[ N^{\beta}, 8 N^{\beta} \right]} \mathbbm{1} \left( \abs{z_m - z_n} < \frac{16}{2^u} \right) \ll$$
$$\ll N^{\varepsilon} E^{*} (X, N) \ll N^{\theta + \varepsilon}$$
Over all we have bounded our sum by
$$2^u N^{\theta + \varepsilon} \ll N^{4 - \varepsilon}$$
for $2^u \leq N^{4 - \theta - \varepsilon}$, as required.

\section{Tightness of Bounds}
%beta to theta%
Before we prove Theorem \ref{zeta v2} let us show that the bounds we have just shown are optimal, assuming of course the truth of the Lindel\"{o}f Hypothesis. Take the sequence $x_n = \floor{n^{4 - \theta}}$ for some $2 \leq \theta \leq 3$, by Theorem 2 of \cite{Robert&Sargos} we have that the number of solutions to the Diophantine inequality
$$\abs{m_{1}^{4 - \theta} + m_{2}^{4 - \theta} - m_{3}^{4 - \theta} - m_{4}^{4 - \theta}} < \delta N^{4 - \theta}$$
where $1 \leq m_1, m_2, m_3, m_4 \leq N$ is at most $N^{2 + \varepsilon} + \frac{N^{4 + \varepsilon}}{\delta}$. Choosing $\delta = \frac{1}{N^{4 - \theta}}$ we see that in our case we have
$$E^{*} (\floor{n^{4 - \theta}}, N) \ll N^{\theta + \varepsilon}$$
and so by Theorem \ref{zeta v1}, under the Lindel\"{o}f Hypothesis we have for almost $\alpha$
$$\delta_{\min}^{\alpha} (N) \ll \frac{1}{N^{4 - \theta - \varepsilon}}$$
On the other hand, clearly we have
$$\delta_{\min}^{\alpha} (N) \geq \min_{1 \leq n \leq N^{4 - \theta}} \norm{\alpha n}$$
and it is well known that for almost all $\alpha$ we have $\min_{1 \leq n \leq M} \norm{\alpha n} \gg \frac{1}{M^{1 + \varepsilon}}$, therefore for almost all $\alpha$ we have
$$\delta_{\min}^{\alpha} (N) \gg \frac{1}{N^{4 - \theta + \varepsilon}}$$
which shows that we cannot improve the exponent.

\section{Unconditional Results}
In order to show that the minimal gap is $\ll \frac{1}{N^{\eta}}$ for some $\eta$, it is sufficient to show that for $2^u \ll N^{\eta - \varepsilon}$ we have
$$\sum_{2^{u - 1} \leq j_1, j_2 < 2^u} \sum_{z_m, z_n} \mathbbm{1} \left( \abs{j_1 z_m - j_2 z_n} \right) \ll N^{4 - \varepsilon}$$
As in the previous section, we can split the sum over $z_m, z_n$ into a sum over $\beta$ of sums of the form
$$\sum_{2^{u - 1} \leq j_1, j_2 < 2^u} \sum_{z_m, z_n \in \left[ N^{\beta}, 8 N^{\beta} \right]} \mathbbm{1} \left( \abs{j_1 z_m - j_2 z_n} \right)$$
However, since now we are no longer assuming that our sequence is of polynomial growth, there may be many intervals. To overcome this, let $L \leq N^3$ be some number we will choose later. We will split the range of $z_m, z_n$ into $\log L \ll \log N$ intervals of the form $z_m, z_n \in \left[ N^{\beta}, 8 N^{\beta} \right]$ where $N^{\beta} \leq L$ and into $z_m, z_n \geq L$, and we will show that the sum in each of these is $\ll N^{4 - \varepsilon}$.

For $z_m, z_n \in \left[ N^{\beta}, 8 N^{\beta} \right]$ we can use the same setup as in the previous section. The argument follows through up to our bound on
$$\mathrm{I_1} = \int_{t \geq \frac{1}{T}} \int_{\mathbb{R}} \zeta \left( \frac{1}{2} + i t T + i u \right) \zeta \left( \frac{1}{2} - i t T + i u \right) F(u) \ \mathrm{d} u \abs{P(t T)}^2 \Phi(t) \ \mathrm{d} t \ll$$
$$\ll \frac{1}{T} \int_{\mathbb{R}} \frac{1}{u^2 + 1} \int_{\mathbb{R}} \abs{\zeta \left( \frac{1}{2} + i t + i u \right)} \abs{\zeta \left( \frac{1}{2} - i t + i u \right)} \abs{P(t)}^2 \Phi \left( \frac{t}{T} \right) \mathrm{d} t \ \mathrm{d} u$$
By the convexity bound $\zeta \left( \frac{1}{2} + i t \right) \ll \abs{t}^{\frac{1}{4} + \varepsilon}$ we see that the contribution of $u \geq t$ is negligible. For the rest of the integral, using H\"{o}lder's inequality with $\frac{2}{\kappa} + \frac{1}{B} = 1$ we get
$$\int_{\mathbb{R}} \abs{\zeta \left( \frac{1}{2} + i t + i u \right)} \abs{\zeta \left( \frac{1}{2} - i t + i u \right)} \abs{P(t)}^2 \Phi \left( \frac{t}{T} \right) \ll$$
$$\ll \left( \int_{\mathbb{R}} \abs{\zeta \left( \frac{1}{2} + i t + i u \right)}^{\kappa} \Phi \left( \frac{t}{T} \right) \mathrm{d} t \right)^{\frac{1}{\kappa}} \left( \int_{\mathbb{R}} \abs{\zeta \left( \frac{1}{2} - i t + i u \right)}^{\kappa} \Phi \left( \frac{t}{T} \right) \mathrm{d} t \right)^{\frac{1}{\kappa}} \cross$$
$$\cross \left( \int_{\mathbb{R}} \abs{P(t)}^{2(B - 1)} \abs{P(t)}^2 \Phi \left( \frac{t}{T} \right) \right)^{\frac{1}{B}}$$
By the rapid decay of $\Phi$ the contribution of $t \geq T^{1 + \varepsilon}$ in the first two integrals is negligible, and so we have
$$\left( \int_{\mathbb{R}} \abs{\zeta \left( \frac{1}{2} + i t + i u \right)}^{\kappa} \Phi \left( \frac{t}{T} \right) \mathrm{d} t \right)^{\frac{1}{\kappa}} \left( \int_{\mathbb{R}} \abs{\zeta \left( \frac{1}{2} - i t + i u \right)}^{\kappa} \Phi \left( \frac{t}{T} \right) \mathrm{d} t \right)^{\frac{1}{\kappa}} \ll T^{\frac{2 + 2 m(\kappa)}{\kappa}}$$
As for the third integral, we have
$$\left( \int_{\mathbb{R}} \abs{P(t)}^{2 \left( B - 1 \right)} \abs{P(t)}^2 \Phi \left( \frac{t}{T} \right) \mathrm{d} t \right)^{\frac{1}{B}} \leq \abs{P(0)}^{2 \left( 1 - \frac{1}{B} \right)} \left( \int_{\mathbb{R}} \abs{P(t)}^2 \Phi \left( \frac{t}{T} \right) \mathrm{d} t \right)^{\frac{1}{B}} \ll$$
$$\ll \left( N^{\beta} E^{*} (X, N) \right)^{1 - \frac{1}{B}} \left( T E^{*} (X, N) \right)^{\frac{1}{B}} = N^{\frac{2 \beta}{\kappa}} T^{1 - \frac{2}{\kappa}} E^{*} (X, N)$$
Therefore, we have shown
$$\mathrm{I_1} \ll T^{\frac{2 m(\kappa)}{\kappa}} N^{\frac{2 \beta}{\kappa}} E^{*} (X, N) \leq 2^{\frac{2 u m(\kappa)}{\kappa}} N^{\beta \frac{2 + 2 m(\kappa)}{\kappa}} E^{*} (X, N) \ll 2^{\frac{2 u m(\kappa)}{\kappa}} L^{\frac{2 + 2 m(\kappa)}{\kappa}} E^{*} (X, N)$$
which shows that
$$\sum_{2^{u - 1} \leq j_1, j_2 < 2^u} \sum_{z_m, z_n \leq L} \mathbbm{1} \left( \abs{j_1 z_m - j_2 z_n} < 1 \right) \ll 2^u 2^{\frac{2 u m(\kappa)}{\kappa}} L^{\frac{2 + 2 m(\kappa)}{\kappa}} E^{*} (X, N)$$
Now we need to bound the contribution of $z_m, z_n \geq L$. Let $T = 2^{u - 1} L$. For any integer $k \geq 1$, set $b_k$ to be the number of $z_m$'s in the interval $[k, k + 1)$, that is
$$b_k = \sum_{m} \mathbbm{1} \left( z_m \in [k, k + 1) \right)$$
We have $\sum_{k = 1}^{\infty} b_k = \frac{N ( N - 1)}{2}$. Notice that
$$\sum_{k = 1}^{\infty} b_{k}^2 \ll \sum_{m, n} \mathbbm{1} \left( \abs{z_m - z_n} < 1 \right) = E^{*} (X, N)$$
Split the interval $[1, \infty)$ into disjoint intervals
$$I_h = \left[ e^{\frac{h}{T}}, e^{\frac{h + 1}{T}} \right)$$
for each $h \geq 0$ and set
$$a_h = \left( \sum_{k \in I_h} b_{k}^2 \right)^{\frac{1}{2}}$$
Fix $j_1, j_2$ and assume without loss of generality that $j_2 \geq j_1$. Assume that $\abs{j_1 z_m - j_2 z_n} < 1$ where $z_m \in I_{h_1}, \ z_n \in I_{h_2}$ for some $h_1, h_2 \geq 0$. This implies that
$$\abs{\frac{j_2}{j_1} - \frac{z_m}{z_n}} < \frac{1}{j_1 z_n} \leq \frac{1}{2^{u - 1} L} \leq \frac{1}{T}$$
As before, we have
$$\frac{1}{2} \leq \frac{z_m}{z_n} \leq 2$$
We also clearly have
$$e^{\left( h_1 - h_2 - 1 \right) / T} \leq \frac{z_m}{z_n} \leq e^{\left( h_1 - h_2 + 1 \right) / T}$$
Since $\frac{z_m}{z_n} \leq 2$, then $e^{\left( h_1 - h_2 \right) / T} \leq 2 e^{\frac{1}{T}}$ and so
$$\abs{\frac{z_m}{z_n} - e^{\left( h_1 - h_2 \right) / T}} \leq e^{\left( h_1 - h_2 \right) / T} \left( e^{\frac{1}{T}} - 1 \right) \leq \frac{4}{T}$$
By the triangle inequality we get
$$\abs{\frac{j_2}{j_1} - e^{\left( h_1 - h_2 \right) / T}} \leq \frac{5}{T}$$
which gives us some condition that $h_1, h_2$ must satisfy in order to contribute to the sum. Now fix $h_1, h_2$, and take some $z_n \in [k, k + 1)$ where $k \in I_{h_2}$. In order for $\abs{j_1 z_m - j_2 z_n} < 1$ to hold we must have $\abs{\ceil{\frac{j_2 k}{j_1}} - z_m} < 4$, because $1 \leq \frac{j_2}{j_1} < 2$. As $\frac{j_2}{j_1} \geq 1$, the mapping $k \mapsto \ceil{\frac{j_2 k}{j_1}}$ is injective, therefore
$$\sum_{z_m \in I_{h_1}, z_n \in I_{h_2}} \mathbbm{1} \left( \abs{j_1 z_m - j_2 z_n} < 1 \right) \ll \sum_{k \in I_{h_2}} b_k \sum_{\substack{z_m \in I_{h_1} \\ \abs{\ceil{\frac{j_2 k}{j_1}} - z_n} < 4}} 1 \ll$$
$$\ll \sum_{k \in I_{h_2}} \sum_{\substack{\ell = - 4 \\ \ceil{\frac{j_2 k}{j_1}} + \ell \in I_{h_1}}} b_k b_{\ceil{\frac{j_2 k}{j_1}} + \ell} \ll \left( \sum_{k \in I_{h_2}} b_{k}^2 \right)^{\frac{1}{2}} \left( \sum_{k \in I_{h_1}} b_{k}^2 \right)^{\frac{1}{2}} = a_{h_1} a_{h_2}$$
by Cauchy-Schwarz. Combining what we know, we get that for fixed $j_1, j_2$ we have
$$\sum_{z_m, z_n \geq L} \mathbbm{1} \left( \abs{j_1 z_m - j_2 z_n} < 1 \right) \ll \sum_{\substack{h_1, h_2 \geq 0 \\ \abs{e^{\left( h_1 - h_2 \right) / T} - \frac{j_2}{j_1}} < \frac{5}{T}}} a_{h_1} a_{h_2}$$
and so summing over all $j_1, j_2 \in \left[ 2^{u - 1}, 2^u \right)$ we have that
$$\sum_{2^{u - 1} \leq j_1, j_2 < 2^u} \sum_{z_m, z_n \geq L} \mathbbm{1} \left( \abs{j_1 z_m - j_2 z_n} < 1 \right) \ll$$
$$\ll \sum_{2^{u - 1} \leq j_1, j_2 < 2^u} \sum_{\substack{h_1, h_2 \geq 0 \\ \abs{e^{\left( h_1 - h_2 \right) / T} - \frac{j_2}{j_1}} < \frac{5}{T}}} a_{h_1} a_{h_2}$$
Now we take as before
$$D(t) = \sum_{2^{u - 1} \leq j < 2^u} j^{i t}$$
but now
$$P(t) = \sum_{h = 0}^{\infty} a_h e^{i h t / T}$$
The key point here is that as before, the constant coefficient of $\abs{P(t)}^2$ counts something very similar to the additive energy, but now all non-constant exponents have a phase of absolute value at least $\frac{1}{T}$ and so it is easier to precisely pick out the constant coefficient. By a similar reasoning to our argument above, letting $\Phi$ be a Schwarz function whose Fourier transform is non-negative, compactly supported on $[-7, 7]$ and bounded there by $2$, and for all $x \in [-6, 6]$ we have $\hat{\Phi} (x) \geq 1$ we have
$$\sum_{2^{u - 1} \leq j_1, j_2 < 2^u} \sum_{\substack{h_1, h_2 \geq 0 \\ \abs{e^{\left( h_1 - h_2 \right) / T} - \frac{j_2}{j_1}} < \frac{5}{T}}} a_{h_1} a_{h_2} \leq \int_{\mathbb{R}} \abs{D(t T)}^2 \abs{P(t T)}^2 \Phi (t) \mathrm{d} t$$
We can do the analysis exactly as above, except for two differences: the first difference is when we pointwise bound for $\abs{t} \leq \frac{1}{T}$ we only have the trivial bound
$$\abs{P(t T)}^2 \leq \abs{P(0)}^2 \leq \left( \sum_{h = 0}^{\infty} a_h \right)^2 \leq \left( \sum_{k = 1}^{\infty} b_k \right)^2 \leq N^4$$
and so the contribution of the integral for $\abs{t} \leq \frac{1}{T}$ is bounded by
$$\frac{2^{2 u} N^4}{T} = \frac{2^{2 u} N^4}{2^{u - 1} L}$$
We see that in order for this to be $\ll N^{4 - \varepsilon}$ we must have $L \gg 2^u N^{\varepsilon}$. The second difference is that there is an extra step when we compute the $L^2$ norm of $\abs{P(t)}^2 \Phi \left( \frac{t}{T} \right)$:
$$\int_{\mathbb{R}} \abs{P(t)}^2 \Phi \left( \frac{t}{T} \right) \mathrm{d} t \ll T \sum_{\substack{h_1, h_2 \\ \abs{h_1 - h_2} \leq 7}} a_{h_1} a_{h_2} \ll T \sum_{h = 0}^{\infty} a_{h}^2 = T \sum_{k = 1}^{\infty} b_{k}^2 \ll T E^{*} (X, N)$$
Over all, we get that
$$\sum_{2^{u - 1} \leq j_1, j_2 < 2^u} \sum_{z_m, z_n \geq L} \mathbbm{1} \left( \abs{j_1 z_m - j_2 z_n} < 1 \right) \ll 2^u 2^{\frac{2 u m(\kappa)}{\kappa}} L^{\frac{2 m(\kappa)}{\kappa}} N^{\frac{8}{\kappa}} E^{*} (X, N)^{1 - \frac{2}{\kappa}}$$
provided that $L \gg 2^u N^{\varepsilon}$. We see that our bounds are better the smaller $L$ is, and so choosing $L = 2^u N^{\varepsilon}$ we get
$$\sum_{2^{u - 1} \leq j_1, j_2 < 2^u} \sum_{m, n} \mathbbm{1} \left( \abs{j_1 z_m - j_2 z_n} < 1 \right) \ll \left( 2^u \right)^{1 + \frac{4 m(\kappa)}{\kappa}} N^{\varepsilon} \left( 2^{\frac{2 u}{\kappa}} E^{*} (X, N) + N^{\frac{8}{\kappa}} E^{*} (X, N)^{1 - \frac{2}{\kappa}} \right)$$
In the case where $E^{*} (X, N) \ll N^{\theta + \varepsilon}$ an easy computation shows that this is $\ll N^{4 - \varepsilon}$ if $2^u \ll N^{\frac{\left( 4 - \theta \right) \left( \kappa - 2 \right)}{2 \kappa - 4} - \varepsilon}$, as required.

Let us just mention that in principle we could have done the whole argument with $\sum_{h = 0}^{\infty} a_h e^{i h t / T}$ instead of $\sum_{n} z_{n}^{i t}$, but we chose not to do so in order to make things more transparent.

\section{Rapidly Growing Sequences}
In this section we will prove Theorem \ref{rapid}. We assume that $x_{n+1} > x_n \left( 1 + \frac{1}{n^{o(1)}} \right)$. Fix $\delta > 0$. We will show that there is some small constant $\varepsilon > 0$ such that for every $u$ that satisfies $2^{u} < N^{2 - \delta}$, we have
$$\sum_{2^{u - 1} \leq j_1, j_2 < 2^u} \sum_{m, n} \mathbbm{1} \left( \abs{j_1 z_m - j_2 z_n} < 1 \right) \ll N^{4 - \varepsilon}$$
which as we have shown is sufficient to deduce that for almost all $\alpha$ we have $\delta_{\min}^{\alpha} (N) < \frac{1}{N^{2 - o(1)}}$.

Call a pair $(m, n)$ \textit{admissible} if there exist $2^{u - 1} \leq j_1, j_2 < 2^u$ such that $\abs{j_1 z_m - j_2 z_n} < 1$, or equivalently $\abs{\frac{z_m}{z_n} - \frac{j_2}{j_1}} < \frac{1}{j_1 z_n}$. Recall that $j_2, j_1$ are in a dyadic interval, and so $\frac{j_2}{j_1} \in \left[ \frac{1}{2}, 2 \right]$, and in particular any admissible pair satisfies $\frac{z_m}{z_n} \in \left[ \frac{1}{3}, 3 \right]$.

By our growth conditions, for any constant $C > 0$ there are at most $N^{\frac{\delta}{2}} \ n$'s such that $z_n < N^C$, and so there are at most $N^{\delta}$ admissible pairs $(m, n)$ such that $z_m < N^C$ or $z_n < N^C$. Each admissible pair contributes at most $2^{2 u} < N^{4 - 2 \delta}$ to the sum (as that is the number of pairs $j_1, j_2$ in the allowed range), and so the contribution of $m, n$ such that $z_m < N^C$ or $z_n < N^C$ is at most $N^{4 - \delta}$.

From here on we assume that $z_m, z_n > N^C$ for some arbitrarily large constant $C$. Take an admissible pair $(m, n)$, and take two pairs $(j_1, j_2), (j_{1}^{'}, j_{2}^{'})$ which satisfy the relevant inequality, that is
$$\abs{\frac{z_m}{z_n} - \frac{j_2}{j_1}} < \frac{1}{j_1 z_n}$$
$$\abs{\frac{z_m}{z_n} - \frac{j_{2}^{'}}{j_{1}^{'}}} < \frac{1}{j_{1}^{'} z_n}$$
By the triangle inequality,
$$\abs{\frac{j_2}{j_1} - \frac{j_{2}^{'}}{j_{1}^{'}}} < \frac{1}{z_n}$$
However, the LHS is a fraction with denominator $< N^{4 - 2 \delta}$, and as $z_n > N^C$ we must have
$$\frac{j_2}{j_1} = \frac{j_{2}^{'}}{j_{1}^{'}}$$.

In particular, this implies that each admissible pair $(m, n)$ contributes at most $N^{2 - \delta}$ to our sum, and so if we show that there are at most $N^{2 + \frac{\delta}{2}}$ admissible pairs our sum will be of size at most $N^{4 - \frac{\delta}{2}}$, which will prove our claim.

Fix $z_n = x_{n_1} - x_{n_2}$. We will show that there are at most $N^{\frac{\delta}{2}} \ m$'s such that the pair $(m, n)$ is admissible, and as there are $< N^2$ options for $n$, this will prove what we need. Take $m$ such that $(m, n)$ is admissible, and write $z_m = x_{m_1} - x_{m_2}$. By our growth conditions, for some arbitrarily small $\varepsilon$ we have
$$\frac{x_{n_1}}{n_{1}^{\varepsilon}} \ll z_n < x_{n_1}$$
$$\frac{x_{m_1}}{m_{1}^{\varepsilon}} \ll z_m < x_{m_1}$$
And so remembering that
$$1 \gg \frac{z_m}{z_n} \gg 1$$
we get that
$$\frac{x_{n_1}}{n_{1}^{\varepsilon}} \ll x_{m_1} \ll n_{1}^{\varepsilon} \cdot x_{n_1}$$
By the fact that $\left( 1 + \frac{1}{n^\varepsilon} \right)^{n^{\varepsilon} \log (n)} \gg n > n^{\varepsilon}$ we see easily that $\abs{m_1 - n_1} \ll n_{1}^{3 \varepsilon} < N^{2 \varepsilon}$, and so there are at most $N^{2 \varepsilon}$ options for $m_1$.
Now fix $m_1$, and we shall show that there are at most $N^{2 \varepsilon}$ options for $m_2$. Assume by contradiction that this is not the case. Then there exist $m_2 < m_{2}^{'} < m_1 - N^{2 \varepsilon}$ such that $(m, n), (m', n)$ are both admissible pairs, where
$$z_m = x_{m_1} - x_{m_2}, \ z_{m'} = x_{m_1} - x_{m_{2}^{'}}$$
That is there exist $2^{u - 1} \leq j_1, j_2, k_1, k_2 < 2^u$ such that
$$\abs{\frac{z_m}{z_n} - \frac{j_2}{j_1}} < \frac{1}{j_1 z_n}$$
$$\abs{\frac{z_{m'}}{z_n} - \frac{k_2}{k_1}} < \frac{1}{k_1 z_n}$$
By using the triangle inequality, we see that
$$\abs{\left( \frac{k_2}{k_1} - \frac{j_2}{j_1} \right) - \frac{z_{m'} - z_m}{z_n}} < \frac{1}{z_n}$$
Now using the reverse triangle inequality $\abs{z - w} \geq \abs{z} - \abs{w}$ we get that
$$\abs{\frac{k_2}{k_1} - \frac{j_2}{j_1}} < \frac{x_{m_{2}^{'}} - x_{m_2}}{z_n} + \frac{1}{z_n}$$
As we have shown,
$$x_{m_1} \ll m_{1}^{\varepsilon} \cdot z_n$$
and so
$$x_{m_2} \ll \frac{x_{m_1}}{\left( 1 + \frac{1}{m_{1}^{\varepsilon}} \right)^{N^{2 \varepsilon}}} < \frac{z_n}{2 N^4}$$
Combining this information we get
$$\abs{\frac{k_2}{k_1} - \frac{j_2}{j_1}} < \frac{1}{N^4}$$
However the LHS is a fraction with denominator smaller than $N^4$, which means that it is $0$. Using our triangle inequality from the beginning again we get
$$\frac{x_{m_{2}^{'}} - x_{m_2}}{z_n} < \frac{1}{z_n}$$
and so
$$x_{m_{2}^{'}} - x_{m_2} < 1$$
a clear contradiction to our growth conditions.
All in all, we have seen that for each $n$ there are at most $N^{3 \varepsilon} \ m$'s such that $(n, m)$ is an admissible pair. Choosing $\varepsilon = \frac{\delta}{6}$ we get our desired bound. $\blacksquare$

\section{Quadratic Polynomials}
In this section we will prove Theorem \ref{quadratic}. The gaps of quadratic polynomials are well studied, see for example \cite{Rudnick-Sarnak, Zaharescu}. In the case where $a, b, c$ are integers (or more generally, rationals) this follows from \cite{Integer Minimal Gaps}, and in fact Theorem 1.2 in \cite{Rectangular Billiard} shows a stronger statement: for almost all $\alpha$, not only do there exist $1 \leq m \neq n \leq N$ and some integer $k$ such that $\abs{\alpha(m^2 - n^2) - k} < 1$, but we can choose $k = k_1 k_2$ where $c(\alpha) N \leq k_1, k_2 \leq C(\alpha) N$. However, it appears that for general real numbers $a, b, c$ this is a new result.

First of all, notice that we can assume without loss of generality that $a = 1, c = 0$ by rescaling and translating the sequence, and so we can work with the sequence $n^2 + \beta n$ for some $\beta \in \mathbb{R}$. The gaps of this sequence are of the form $n^2 + \beta n - m^2 - \beta m = (n - m)(n + m + \beta)$ and so it is enough to bound
$$\sum_{1 \leq j_1, j_2 < M} \sum_{z_m, z_n} \mathbbm{1} \left( \abs{j_1 z_m - j_2 z_n} < 1 \right) \ll \sum_{1 \leq j_1, j_2 < M} \sum_{1 \leq x < y, z < w \leq 2 N} \mathbbm{1} \left( \abs{j_1 x(y + \beta) - j_2 z(w + \beta)} < 1 \right)$$
We need to show that this sum is at most $\frac{N^4}{\log (N)^{1 + \varepsilon}}$. Notice that
$$j_1 x (y + \beta) - j_2 z (w + \beta) = j_1 x y - j_2 z w + \beta (j_1 x - j_2 z)$$
Fix $k = j_1 x - j_2 z$. Then there are at most 2 options for $j_1 x y - j_2 z w$, which we will denote by $C$. Then we have
$$j_1 x y - j_2 z w = j_2 z (y - w) + k y = C$$
and so
$$j_2 z (y - w) = C - k y$$
Fixing $y$, there are at most $d(C - k y)^2$ options for $z, y - w, j_2$, but using $w = y - (y - w), \ j_1 x = j_2 z + k$ we see that we can find $w, j_1 x$ as well. Now, there are $d(j_2 z + k)$ options for $j_1, x$ and so over all using the standard divisor bound $d(n) \ll n^{\varepsilon}$ we have shown that fixing $k, y$ there are $\ll N^{\varepsilon}$ options for $j_1, j_2, x, y, z, w$. As $k \leq M N, \ y \leq 2 N$, running over all pairs of $k, y$, the sum is at most $M N^{2 + \varepsilon}$, and choosing $M = N^{2 - 2 \varepsilon}$ we get what we desire.

Of course, a slightly stronger result could be achieved by using a stronger divisor bound like
$$d(n) \ll 2^{(1 + \varepsilon) \frac{\log N}{\log \log N}}$$

In the case where our sequence is a quadratic polynomial with integer (or rational) coefficients, we can divide the sum according to the value of $a b$ and use the fact that $\sum_{n \leq N^2} d(n)^2 \ll N^2 \log(N)^3$ to get a minimal gap of size
$$\frac{\log (N)^{4 + \varepsilon}}{N^2}$$
which is optimal up to a power of $\log N$.

\section{Minimal Gap in $(\sqrt{n})_{n \geq 1}$}
In this section we will discuss the fractional parts of the sequence $\sqrt{n}$ and prove Theorem \ref{square root}. Of course, technically the minimal gap is $0$ as $\sqrt{4} - \sqrt{1} = 1$, and so we throw away all of the integers in this sequence, that is square roots of perfect squares. Then, defining
$$\delta_{\min}^{2}(N) = \min \left( \norm{\sqrt{m} - \sqrt{n}}: m \neq n \leq N, \sqrt{m}, \sqrt{n} \not\in \mathbb{Z} \right)$$
We shall show that
$$\delta_{\min}^{2}(N) \sim \frac{1}{2 N^{\frac{3}{2}}}$$
This is in stark contrast to the random case where we expect $\delta_{\min}^{2}(N)$ to be on order of $\frac{1}{N^2}$.

Take $n > m$ achieving the minimum in $\delta_N$. Then there exists $k \in \mathbb{Z}$ and $\varepsilon \in \mathbb{R}$ such that
$$\norm{\sqrt{n} - \sqrt{m}} = \abs{\varepsilon}$$
and
$$\sqrt{n} - \sqrt{m} = k + \varepsilon$$
Squaring both sides we get
$$n + m - \sqrt{4 m n} = k^2 + 2 k \varepsilon + \varepsilon^2$$
or equivalently
$$\sqrt{4 m n} - \left( m + n - k^2 \right) = 2 k \varepsilon + \varepsilon^2$$
Notice that clearly the right hand side is non zero, as $\varepsilon \neq 0, -2 k$ and so the left hand side is non zero as well. Then, we have
$$\abs{\sqrt{4 m n} - \left( m + n - k^2 \right)} \geq \abs{\sqrt{4 m n} - \sqrt{4 m n + 1}} = \frac{1}{\sqrt{4 m n} + \sqrt{4 m n + 1}} = \frac{1 + o(1)}{4 \sqrt{m n}}$$
and so we have
$$2 k \abs{\varepsilon} = (1 + o(1)) \abs{2 k \varepsilon + \varepsilon^2} \geq \frac{1 + o(1)}{4 \sqrt{m} \sqrt{n}}$$
and therefore using the arithmetic mean-geometric mean inequality we get
$$\delta_{\min}^{2}(N) = \abs{\varepsilon} \geq \frac{1 + o(1)}{8 \sqrt{n} \sqrt{m} (\sqrt{n} - \sqrt{m})} \geq \frac{1 + o(1)}{2 n^{\frac{3}{2}}} \geq \frac{1 + o(1)}{2 N^{\frac{3}{2}}}$$
which proves the lower bound.

On the other hand, choose $d$ to be the unique positive integer such that
$$4(d + 1)^2 > N \geq 4 d^2$$
and take
$$n = \sqrt{4 d^2 - 2}, m = \sqrt{d^2 - 1}$$
We will use the Taylor series expansion of $\sqrt{x}$ around $1$ to the second order, which is
$$\sqrt{1 + x} = 1 + \frac{1}{2} x - \frac{1}{8} x^2 + O(x^3)$$
and get that
$$\sqrt{m} = d \sqrt{1 - \frac{1}{d^2}} = d - \frac{1}{2 d} - \frac{1}{8 d^3} + O \left( \frac{1}{d^5} \right)$$
$$\sqrt{n} = 2 d \sqrt{1 - \frac{1}{2 d^2}} = 2d - \frac{1}{2 d} - \frac{1}{16 d^3} + O \left( \frac{1}{d^5} \right)$$
and so
$$\sqrt{n} - \sqrt{m} = d + \frac{1}{16 d^3} + O \left( \frac{1}{d^5} \right)$$
and therefore
$$\delta_{\min}^{2}(N) \leq \frac{1 + o(1)}{16 d^3} = \frac{1 + o(1)}{2 N^{\frac{3}{2}}}$$
which is our required lower bound. Combining these, we have shown that
$$\delta_{\min}^{2}(N) \sim \frac{1 + o(1)}{2 N^{\frac{3}{2}}}$$
as requested.

Note that this exact argument shows in fact that if $a, b, c \leq N$ then
$$\min_{\substack{1 \leq a, b, c \leq N \\ \sqrt{a} - \sqrt{b} - \sqrt{c} \not\in \mathbb{Z}}} \left( \abs{\sqrt{a} - \sqrt{b} - \sqrt{c}} \right) \sim \frac{1}{2 N^{\frac{3}{2}}}$$
which is a strictly stronger statement ($\ref{square root}$ is this statement when we only let $c$ range over squares).

One can also ask about the minimal gap in $n^{\frac{1}{k}}$ for $1 < k \in \mathbb{Z}$, that is looking at
$$\delta_{\min}^{k}(N) = \min \left( \norm{\sqrt[k]{n} - \sqrt[k]{m}}: m \neq n \leq N, \sqrt[k]{n}, \sqrt[k]{m} \not\in \mathbb{Z} \right)$$
The case $k = 2$ is the case we have just considered. Our example easily generalizes via
$$n = 2^k d^k - 2^{k-1}, m = d^k - 1$$
to show that
$$\delta_{\min}^{k}(N) \ll \frac{1}{N^{2 - \frac{1}{k}}}$$
However our method to achieve a lower bound generalizes as follows:
$$1 \leq N_{\mathbb{Q} \left( \sqrt[k]{m}, \sqrt[k]{n} \right) / \mathbb{Q}} \left( \sqrt[k]{m} - \sqrt[k]{n} - \ell \right) \ll \left( \sqrt[k]{m} - \sqrt[k]{n} - \ell \right) N^{k - \frac{1}{k}}$$
which gives just
$$\delta_{\min}^{k}(N) \gg \frac{1}{N^{k - \frac{1}{k}}}$$
which for $k > 2$ does not match our lower bound. It would be interesting to know the answer for general $k$.

\section{Sequences with large additive energy}
In this section we will discuss sequences with large additive energy. For simplicity, we will assume that our sequences are integer-valued, though under some strong assumptions we can work with real-valued sequences as well. For a set $A$ of integers we recall that its additive energy is defined by
$$E(A) = \# \left\{ (a, b, c, d) \in A^4 \ \big\rvert \ a + b = c + d \right\}$$

In \cite{Integer Minimal Gaps} it is shown that a sequence with additive energy of $\ll N^{2 + \varepsilon}$ (the lowest possible energy) has a metric minimal gap of size at most $\frac{1}{N^{2 - \varepsilon}}$. One can ask the opposite question: does a large additive energy imply that the sequence has a large minimal gap? The answer to this is easily seen to be no. For example, take the sequence defined by
$$x_{2 k} = k^2, \ x_{2 k + 1} = k$$
Then this sequence has a small metric minimal gap, because $k^2$ has a small metric minimal gap, but its additive energy is of order $N^3$.

Let us look more carefully at our argument for Theorem \ref{lower bound}. Denoting the set $X_N = \{ x_1, \dots, x_N \}$ We used the random variable
$$\sum_{m, n \in X_N, m \neq n} \mathbbm{1}_{[- \frac{1}{M}, \frac{1}{M}]} \left( \alpha x_m - \alpha x_n \right)$$
However, there may be multiplicities in the values of $x_m - x_n$, and then we are summing over some indicators multiple times, and so it would be more economical to work with the random variable
$$\sum_{0 \neq t \in X_N - X_N} \mathbbm{1}_{[- \frac{1}{M}, \frac{1}{M}]} \left( \alpha t \right)$$
Then we see that the expectation of this variable is
$$\frac{2 \abs{X_N - X_N}}{M}$$
and the same type of argument as in section \ref{Lower Bound} shows that in this case the minimal gap is at least of order
$$\frac{1}{\abs{X_N - X_N} \log (N)^{1 + \varepsilon}}$$
Therefore, if $\abs{X_N - X_N}$ is small then the minimal gap is large. By Cauchy-Schwarz one can see that $E(X_N) \geq \frac{\abs{X_N}^4}{\abs{X_N - X_N}}$, and so $\abs{X_N - X_N}$ small implies that $E(X_N)$ is large. We can use the following theorems from additive combinatorics to give a partial converse to this.

\begin{thm}\label{Balog}

(Balog-Szemer\'{e}di-Gowers Theorem)

Let $A$ be a finite subset of an abelian group. If for some $K$ we have $E(A) \geq \frac{A^3}{K}$ (where $E(A)$ is the additive energy) then there is a subset $A' \subset A$ such that $\abs{A'} \geq K^{- \mathcal{O}(1)} \abs{A}$ and $\abs{A' + A'} \leq K^{\mathcal{O}(1)} \abs{A'}$

\end{thm}

\begin{thm}\label{Plun}

(Plünnecke-Ruzsa inequality)

If $A, B$ are finite subsets of an abelian group and $K$ is a constant such that $\abs{A + B} \leq K \abs{A}$, then for all non-negative integers $m, n$ we have $\abs{m B - n B} \leq K^{m + n} \abs{A}$

\end{thm}

Using $A = B = A', m = n = 1$ in \ref{Plun} together with \ref{Balog} we see that if $E(A) \geq \frac{A^3}{K}$ there exists a subset $A' \subset A$ such that $\abs{A'} \geq K^{- \mathcal{O}(1)} \abs{A}$ and $\abs{A' - A'} \leq K^{\mathcal{O}(1)} \abs{A'}$. With this, we can prove Theorem \ref{Large Energy} which we restate here for convenience:

\begin{thm}\label{reLarge Energy}

Let $x_n$ be a monotone increasing sequence of positive integers such that for some $r > 0$
$$E_N := \# \left\{ (a, b, c, d) \in [1, N]^4 \ \big\rvert \ x_a + x_b = x_c + x_d \right\} \geq r N^3$$
Then there exist increasing sequences of positive integers $n_k, a_k \to \infty$ such that $n_k$ is of density bounded from below in $\left\{ 1, \dots, a_N \right\}$ and for almost all $\alpha$ we have
$$\delta_{\min}^{\alpha} (x_{n_k}, a_N) \gg \frac{1}{a_N \log (a_N)^{1 + \varepsilon}}$$

\end{thm}

Using the Cantelli argument as above it is enough to show that for the aforementioned $a_k, x_{n_k}$ we have a small difference set, that is
$$\# \left\{ x_{n_a} - x_{n_b} \ \big\rvert \ 1 \leq n_a, n_b \leq a_N \right\} \leq C(r) a_N \log(a_N)^{\varepsilon}$$
Let $a_0 = 0, a_{n + 1} = 2^{2^{a_n}}$, and divide the positive integers into intervals $I_h = \left( a_h, a_{h + 1} \right]$ for $h \geq 0$. As we have stated, there exists a subset $J_h \subset I_h$ of density at least $c(r) > 0$ such that $A_h = \{ x_k \ \big\rvert \ k \in J_h \}$ satisfies $\abs{A_h - A_h} \leq C(r) \abs{I_h}$. The desired subsequence $n_k$ be the union of all the $J_h$.

Splitting the set of differences into differences where both elements are contained in $I_{N - 1}$ (which are the majority) and differences where this is not the case we get that
$$\# \left\{ x_{n_a} - x_{n_b} \ \big\rvert \ 1 \leq n_a, n_b, \leq a_N \right\} \leq \# \left\{ x_{n_a} - x_{n_b} \ \big\rvert \ a_{N - 1} < n_a, n_b \leq a_N \right\} + 2 a_{N - 1} a_N \ll$$
$$\ll C(r) a_N + a_N \log \log(a_N) \ll a_N \log(a_N)^{\varepsilon}$$
as required.

Of course, under weaker assumptions on the additive energy in Theorem \ref{reLarge Energy} we can get a correspondingly weaker result regarding the subsequence, using exactly the same techniques.

A natural question is how to extend this result to real-valued sequences. Firstly, an important concept will be the covering number of a set. Given a set $A \subset \mathbb{R}$ and a scale $\gamma > 0$, define the covering number of a set to be
$$\mathrm{cov} \left( A, \gamma \right)= \inf_{B \subset \mathbb{R}} \left\{ \abs{B} \ \big\rvert \ \forall a \in A \ \exists b \in B: \abs{a - b} < \gamma \right\}$$
that is, the smallest number of balls with radius $\gamma$ that are needed to cover $A$. This number essentially measures the size of the set $\abs{A}$ if we say that two numbers which are of difference $\ll \gamma$ are equivalent, so this is a "fuzzy" notion of size, and $\gamma$ is the scale of fuzziness which we permit. We can define the additive energy of a set with respect to a certain scale $\gamma$ by
$$E^{*} (A, \gamma) = \# \{ (a, b, c, d) \in A^4 \ \big\rvert \ \abs{a + b - c - d} < \gamma \}$$ Replacing each element $a \in A$ with $\gamma \floor{\frac{a}{\gamma}}$ and then using the above argument we see that given a set $A$ such that $E^{*} (A, \gamma) \gg N^3$ there is a subset $A' \subset A$ such that $\mathrm{cov} \left( A' - A', 2 \gamma \right) \leq K^{\mathcal{O} (1)} \abs{A}$. Using the same argument we did to prove Theorem \ref{reLarge Energy} we see that as before given a sequence $x_n$ of real numbers with large additive energy of scale $\gamma$, there is a subsequence $x_{n_k}$ such that along a certain subsequence $N_k \to \infty$ the difference set has a small size of scale $2 \gamma$.

Now, we can ask the following question: given a sequence $x_n$ whose difference set is of small size (that is, $\mathrm{cov} \left( X_N - X_N, \gamma \right) \ll N$), can we say that the metric minimal gap of $x_n$ is large? This is exactly the content of Theorems \ref{Fine Scale} and \ref{Coarse Scale}, which we will now prove, and once again restate them for convenience:

\begin{thm}\label{reFine Scale}

For any sequence of real numbers $x_n$ such that $x_{n + 1} - x_n \geq \frac{c \log n}{n}$ for some constant $c > 0$, if for some $\varepsilon_0 > 0$ we have
$$\mathrm{cov} \left( X_N - X_N, \gamma \right) \ll N$$
where $\gamma = \frac{1}{N \log(N)^{1 + \varepsilon_0}}$, we have that for every $\varepsilon > 0$ and for almost all $\alpha \in \mathbb{R}$,
$$\delta_{\min}^{\alpha} (x_n) \gg \frac{1}{N \log(N)^{1 + \varepsilon}}$$

\end{thm}

\begin{thm}\label{reCoarse Scale}

There exists a sequence $\varepsilon_n \downarrow 0$ such that the sequence $x_n = n - \frac{1}{n^{1 - \varepsilon_n}}$ satisfies for all $\varepsilon > 0$ and for almost all $\alpha \in \mathbb{R}$ that
$$\delta_{\min}^{\alpha} (N) \ll \frac{1}{N^{2 - \varepsilon}}$$
where the implied constant depends on $\alpha$ and $\varepsilon$.

\end{thm}

Before we prove these, let us briefly give some intuition. Imagine that our sequence was of the form $x_n = a_n + y_n$, where $a_n \in \mathbb{Z}$ and $y_n$ is very small. Then, as $x_n - x_m = a_n - a_m + (y_n - y_m)$, the difference set of $x_n$ would be like the difference set of $a_n$, but with some "fuzz" concentrated around each element of the form $a_n - a_m$, corresponding to the possible values of $y_n - y_m$ when we fix $a_n - a_m$. Theorem \ref{reFine Scale} says that if the $y_n$ are small enough, this "fuzz" is negligible and the answer is the same as in the case where our sequence is integer-valued. To prove Theorem \ref{reCoarse Scale}, we take a sequence of the form $x_n = n + y_n$, such that the difference set is of the form $n - m + \left( y_n - y_m \right)$, and we choose $y_n$

Theorem \ref{reFine Scale} follows from what we have said about integer-valued sequences: take a set $B_N$ of size $\ll N$ such that for every $m, n \leq N$ we have $\abs{x_m - x_n - b} < \frac{1}{N \log(N)^{1 + \varepsilon_0}}$ for some $b \in B$. Then, suppose that
$$\delta_{\min}^{\alpha} (x_n, N) < \frac{1}{N \log(N)^{1 + \varepsilon}}$$
that is, there exist some $m \neq n \leq N$ such that $\norm{\alpha (x_m - x_n)} < \frac{1}{N \log(N)^{1 + \varepsilon}}$, and so $\norm{\alpha b} < \frac{\abs{\alpha} + 1}{N \log(N)^{1 + \varepsilon}}$. Now using the same arguments as in sections 2 and \ref{Lower Bound} shows that this happens for a set of $\alpha$ of measure 0.

Theorem \ref{reFine Scale} combined with our argument above immediately proves Corollary \ref{Corr}.

Now we will prove Theorem \ref{reCoarse Scale}. Obviously for any such sequence where $\varepsilon_n \to 0$, the set of differences is $\frac{1}{N^{1 - \varepsilon}}$ close to the set of differences of $x_n = n$, which is of linear size, and so we just need to show that the minimal gap is sufficiently small. In our sequence we will let $\varepsilon_n = \eta_k$ be constant on intervals of the form $\left( 4^{k - 1}, 4^k \right]$ and $\eta_k$ will decay slowly as $k \to \infty$. Throughout, we will assume that $\eta_k$ is small. We will show that $\delta_{\min}^{\alpha} (N)$ is small for almost all $\alpha \in [-C, C]$ for any $C > 0$, and therefore it is small for almost all $\alpha \in \mathbb{R}$. Notice that $\delta_{\min}^{\alpha} (N) \leq \delta_{\min}^{\alpha} \left( 4^{\floor{\log_4 (N)}} \right)$. Therefore, to show that $\delta_{\min}^{\alpha} (N)$ is small for all $N$ it is enough to look at the case where $N = 4^k$. Clearly,
$$\delta_{\min}^{\alpha}(N) = \min \left\{ \norm{\alpha (x_n - x_m)} \ \big\rvert \ 1 \leq m \neq n \leq N \right\} \leq \min \left\{ \norm{\alpha (x_n - x_m)} \ \big\rvert \ \frac{N}{4} < m \neq n \leq N \right\} =$$
$$\min \left\{ \norm{\alpha \left( n - m - \frac{1}{n^{1 - \eta_k}} + \frac{1}{m^{1 - \eta_k}} \right)} \ \big\rvert \ \frac{N}{4} < m \neq n \leq N \right\}$$
Next we will state a useful lemma:

\begin{lem}\label{Real Coarse Scale}
For all $\eta > 0$ small enough (say $\eta < \frac{1}{10}$), for any $\alpha$ such that for some coprime integers $a, q, \ 0 < q < \frac{N}{3}$ that satisfy
$$\abs{\alpha - \frac{a}{q}} \leq \frac{3}{q N}$$
and
$$\abs{\alpha - \frac{a}{q}} \geq \frac{9}{2 q^{1 + \frac{1}{1 - \eta}}}$$
we have
$$\min \left\{ \norm{\alpha \left( n - m - \frac{1}{n^{1 - \eta}} + \frac{1}{m^{1 - \eta}} \right)} \ \big\rvert \ \frac{N}{4} < m \neq n \leq N \right\} \leq \frac{8}{N^{2 - \eta}}$$
\end{lem}

For now assume this lemma, and we will prove it later. By Dirichlet's approximation theorem, for all $\alpha$ there exist some coprime integers $a, q, \ 0 < q < \frac{N}{3}$ such that
$$\abs{\alpha - \frac{a}{q}} \leq \frac{3}{q N}$$
By Khinchin's theorem, for almost all $\alpha$, for all $q$ large enough we have $\abs{\alpha - \frac{a}{q}} \geq \frac{9}{2 q^2 \log (q)^2}$, and so if $\abs{\alpha - \frac{a}{q}} < \frac{9}{2 q^{1 + \frac{1}{1 - \eta}}}$ we must have
$$q^{\frac{\eta}{1 - \eta}} < \log (q)^2$$
Taking logarithms this is equivalent to
$$\frac{\log q}{\log \log q} < \frac{2 (1 - \eta)}{\eta}$$
and therefore
$$\log q < \frac{4}{\eta} \log \left( \frac{2}{\eta} \right)$$
which means that
$$q < \left( \frac{2}{\eta} \right)^{\frac{4}{\eta}}$$
However, notice that the set of $\alpha \in [-C, C]$ such that $\abs{\alpha - \frac{a}{q}} \leq \frac{1}{q N}$ for some $q < \left( \frac{2}{\eta} \right)^{\frac{4}{\eta}}$ is of measure at most
$$\frac{4 C \left( \frac{2}{\eta} \right)^{\frac{4}{\eta}}}{N}$$
simply by bounding the measure independently for each $q$. That is, for all $\alpha \in [- C, C]$ up to a set of measure at most $\frac{4 C \left( \frac{2}{\eta} \right)^{\frac{4}{\eta}}}{N}$ we have
$$\min \left\{ \norm{\alpha \left( n - m - \frac{1}{n^{1 - \eta}} + \frac{1}{m^{1 - \eta}} \right)} \ \big\rvert \ \frac{N}{4} < m \neq n \leq N \right\} \leq \frac{8}{N^{2 - \eta}}$$

Now we return to our case, where $\eta = \eta_k$. Applying Lemma \ref{Real Coarse Scale} we see that for all $\alpha \in [- C, C]$ up to a set of measure at most $\frac{2 C \left( \frac{2}{\eta_k} \right)^{\frac{4}{\eta_k}}}{4^k}$ we have
$$\delta_{\min}^{\alpha} \left( 4^k \right) \leq \frac{8}{4^{k (2 - \eta_k)}}$$
If we let $\eta_k$ tend to 0 sufficiently slowly so that
$$\sum_{k = 1}^{\infty} \frac{\left( \frac{2}{\eta_k} \right)^{\frac{4}{\eta_k}}}{4^k} < \infty$$
then by Borel-Cantelli, for almost all $\alpha \in [- C, C]$ we will have
$$\delta_{\min}^{\alpha} (4^k) \ll \frac{1}{4^{k (2 - o(1))}}$$
and therefore for almost all $\alpha \in [-C, C]$ we have
$$\delta_{\min}^{\alpha} (N) \ll \frac{1}{N^{2 - o(1)}}$$
as $N \to \infty$, exactly as required. $\blacksquare$
\\ \\
\textbf{Proof of Lemma \ref{Real Coarse Scale}}:

Take some number $\frac{N}{3} \leq k \leq \frac{2 N}{3}$, and let us look at all pairs of $\frac{N}{4} \leq m < n \leq N$ such that $n - m = k$. If $n, m$ is such a pair then obviously so is $n + 1, m + 1$. Notice that
$$\left( x_n - x-m \right) - \left( x_{n + 1} - x_{m + 1} \right) = \left( \frac{1}{m^{1 - \eta}} - \frac{1}{(m + 1)^{1 - \eta}} \right) - \left( \frac{1}{n^{1 - \eta}} - \frac{1}{(n + 1)^{1 - \eta}} \right)$$
This is a positive number which is at most
$$\frac{1}{m^{1 - \eta}} - \frac{1}{(m + 1)^{1 - \eta}} \leq \frac{1 - \eta}{m^{2 - \eta}} < \frac{16}{N^{2 - \eta}}$$
Therefore, the points $x_n - x_m$ for $n - m = k$ fixed are contained in the interval
$$\left[ k + \left( \frac{1}{N - k} \right)^{1 - \eta} - \left( \frac{1}{N} \right)^{1 - \eta}, k + \left( \frac{4}{N} \right)^{1 - \eta} - \left( \frac{1}{k + \frac{N}{4}} \right)^{1 - \eta} \right]$$
Notice that for some $t$ in this interval we have $\alpha t \in \mathbb{Z}$, then there exist some $m, n$ such that $\abs{x_n - x_m - t} < \frac{8}{N^{2 - \eta}}$ and therefore $\norm{\alpha (x_n - x_m)} < \frac{8 \abs{\alpha}}{N^{2 - \eta}}$, which is a small gap. Now we will find such a $t$.

Notice that for $\frac{N}{3} \leq k \leq \frac{2 N}{3}$ (and for $\eta$ close enough to 0) the above interval contains the interval
$$\left[ k + \frac{\left( \frac{3}{2} \right)^{1 - \eta} - 1}{N^{1 - \eta}}, k + \frac{4^{1 - \eta} - \left( \frac{12}{11} \right)^{1 - \eta}}{N^{1 - \eta}} \right]$$
Notice that for $\eta$ small enough this contains the interval
$$\left[ k + \frac{1}{3 N^{1 - \eta}}, k + \frac{7}{3 N^{1 - \eta}} \right]$$
We want $\alpha t$ to be small for some $t$ in one of these intervals, or equivalently $\alpha k \in \left[ - \frac{\alpha(c_1 + c_2)}{N^{1 - \eta}}, - \frac{\alpha c_1}{N^{1 - \eta}} \right]$ for some $k$, where this interval is taken $\bmod \ 1$. Writing $k = \frac{N}{3} + \ell$ for some $0 \leq \ell \leq \frac{N}{3}$ then what we need is
$$\alpha \ell \in \left[ -\frac{N}{3} \alpha - \frac{\alpha}{3 N^{1 - \eta}}, -\frac{N}{3} \alpha - \frac{7 \alpha}{3 N^{1 - \eta}} \right]$$
The exact structure of the right hand side does not matter, what matters is the fact that it is an interval of length $\frac{2 \abs{\alpha}}{N^{1 - \eta}}$ for some absolute constant $c_2$. To emphasize this, let us write that interval as
$$\left[ y, y + \frac{2}{N^{1 - \eta}} \right]$$
for some $y$. Recall that we are given coprime integers $a, q$ such that $0 < q < \frac{N}{3}$ and
$$\abs{\alpha - \frac{a}{q}} \leq \frac{3}{q N}$$
For all $\ell \leq \frac{N}{3}$ we have
$$\abs{\alpha \ell - \frac{a \ell}{q}} \leq \frac{3 \ell}{q N} \leq \frac{1}{q}$$
and so if for some $0 \leq \ell \leq \frac{N}{3}$ we have 
$$\frac{a \ell}{q} \in \left[ y + \frac{1}{q}, y + \frac{2}{N^{1 - \eta}} - \frac{1}{q} \right]$$
then we would finish. Because $a, q$ are coprime, as $\ell$ runs from $0$ to $\frac{N}{3}, \ \frac{a \ell}{q}$ runs over all fractions with denominator $q$. Therefore, if the interval $\left[ y + \frac{1}{q}, y + \frac{2}{N^{1 - \eta}} - \frac{1}{q} \right]$ is of length at least $\frac{1}{q}$ it contains such a fraction which is what we need to prove. Therefore it is sufficient to show that
$$q \geq \frac{3 N^{1 - \eta}}{2}$$
or equivalently
$$\left( \frac{2 q}{3} \right)^{\frac{1}{1 - \eta}} \geq \frac{2 q^{\frac{1}{1 - \eta}}}{3} \geq N$$
Assume to the contrary that this is not the case. Then we have
$$\abs{\alpha - \frac{a}{q}} \leq \frac{3}{q N} < \frac{9}{2 q^{1 + \frac{1}{1 - \eta}}}$$
a contradiction to what we know about $q$. $\blacksquare$.

Notice that Theorem \ref{Fine Scale} is not quite a full converse to \ref{Corr}: in order to prove Corollary \ref{Corr} we used the fact that every sequence with large additive energy has a relatively dense subsequence with a small difference set, and then that every sequence with a small difference set has a large minimal gap. In Theorem \ref{Fine Scale} we showed that the second part of this proof fails in the case where $\gamma = \frac{1}{N^{1 - \varepsilon}}$, as we showed that there exists a sequence with small difference set and large minimal gap. However, theoretically our sequence $n - \frac{1}{n^{1 - \varepsilon_n}}$ could have a subsequence such as described in \ref{Corr}, just we would have to choose the subsequence more carefully. We conjecture that this is not the case. Specifically,

\textbf{Conjecture}: Every subsequence of $n - \frac{1}{n^{1 - \varepsilon_n}}$ with lower positive density has a metric minimal gap of size $\frac{1}{N^{2 - o(1)}}$.
\\ \\
\textbf{Funding}:
\\This research was supported by the European Research Council (ERC) under the European Union’s Horizon 2020 research and innovation programme (Grant agreement No. 786758).
\\ \\
\textbf{Acknowledgements}:
\\This work is part of my M.Sc. thesis. I would like to greatly thank my supervisor, Zeev Rudnick, who was incredibly helpful throughout the writing process. He was always available with useful advice, both mathematical and otherwise, and his many invaluable comments, corrections and suggestions elevated this from a typo-riddled collection of formulae to a hopefully intelligible paper. In short, he taught me how to write an academic paper.

\end{document}